\documentclass[12pt]{article}
\usepackage{hyperref}
\usepackage{float}
\floatstyle{plaintop}
\restylefloat{table}

\usepackage{multirow}
\usepackage[T1]{fontenc}
\usepackage[utf8]{inputenc}
\usepackage{authblk}

\usepackage{dsfont}
\usepackage{amssymb,amsthm,amsmath}
\usepackage{setspace}
\newtheorem{thm}{Theorem}[section]

 \makeatletter

\@addtoreset{lemma}{thm}

\newtheorem{corollary}{Corollary}[section]
 \makeatletter

\@addtoreset{corollary}{thm}

\newtheorem{remark}{Remark}[section]

\numberwithin{equation}{section}
\date{}

\begin{document}

%% The amsthm package provides extended theorem environments

\title{Randomized  pivots for  means  of short and long memory linear processes}

\author{Mikl\'{o}s Cs\"{o}rg\H{o}\thanks{mcsorgo@math.carleton.ca}, Masoud M. Nasari\thanks{mmnasari@math.carleton.ca }  and Mohamedou Ould-Haye\thanks{ouldhaye@math.carleton.ca}  \\
\small{School of Mathematics and Statistics, Carleton University}\\\small{ Ottawa, ON, Canada} }

\maketitle

\begin{abstract}
In this paper we introduce    randomized pivots for the means of short and long memory linear processes. We show that, under the same conditions, these pivots converge in distribution to the same limit as that of their classical non-randomized counterparts. We also present numerical results that indicate that these randomized pivots significantly outperform their classical counterparts and as a result they lead to a more accurate inference about the population mean.
\end{abstract}

\section{Introduction and background}\label{Introduction and Background}
Recently,  Cs\"{o}rg\H{o} and Nasari \cite{Csorgo and Nasari} investigated the problem of establishing   central limit theorems (CLTs) with improved rates  for randomized versions of the  Student $t$-statistic based on  i.i.d. observations $X_1,X_2\ldots$. The resulting  improvements     yield confidence intervals for the population mean $\mu$  with a smaller magnitude of error than that of the classical CLT for the Student $t$-pivot based on   an i.i.d.  sample of size $n$, $n\geq 1$. The improvements  in hand result from incorporating functionals of multinomially distributed random variables as coefficients for, and independent of,  the data. More precisely,    Cs\"{o}rg\H{o} and Nasari \cite{Csorgo and Nasari} introduced and, via conditioning on the random  weights, studied the asymptotic distribution of  the  randomized pivot for the  population mean $\mu:=E X_1$, that is defined  as follows

\begin{equation}\label{added Cs2}
  \frac{\sum_{i=1}^n  \big| \frac{w^{(n)}_i}{m_n}-\frac{1}{n} \big| ( X_i-\mu)   }{S_n \sqrt{\big( \frac{w^{(n)}_i}{m_n}-\frac{1}{n} \big)^2}},
\end{equation}
and  can be computed via generating, independently from the data,   a  realization of the multinomial random weights $(w_{1}^{(n)},\ldots,w_{n}^{(n)})$ with $\sum_{i=1}^n w_{i}^{(n)}=m_n$ and  associated probability vector $(1/n,\ldots,1/n)$; here   $S_{n}^{2}$ is the sample variance.
%We note that one way of generating the random weights $w_i^{(n)}$ independently from the data  is to
%re-sample  with replacement $m_n$ times with respective probabilities $1/n$ from the set of indices $\{1,\ldots,n \}$ of the i.i.d. random variables in the random sample $\{X_1,\ldots,X_n\}$  so that,  for each $1\leq i \leq n$, $w^{(n)}_i$ is the  number of times the index $i$ of $X_i$ in $\{X_i,\ 1\leq i \leq n \}$ is chosen in this re-sampling process.

\par
In paper \cite{Csorgo and Nasari},  it is shown that, on assuming $E|X_1|^3< +\infty$, the magnitude of the error   generated by approximating the sampling distributions of   these   normalized/Studentized  randomized partial sums of i.i.d. observables, as in (\ref{added Cs2}),  by the standard normal distribution function $\Phi(.)$   can be of order $O(1/n)$. The latter rate is achieved when one takes $m_n=n$, and is to be compared to that of the $O(1/\sqrt{n})$ error rate of the Student $t$-statistic under the same moment condition.

\par
%In the present paper the problem of  bootstrapping linear processes is studied.
The present  work is an extension of the results in the  aforementioned paper \cite{Csorgo and Nasari} to short and long memory linear processes via creating  randomized  direct pivots for  the mean  $\mu=E X_1$  of   short and long memory linear processes \emph{\`{a} la} (\ref{added Cs2}). Adaptation of the randomized version of the Student $t$-statistic as in \cite{Csorgo and Nasari}

\begin{equation}\label{added 1}
\frac{\sum_{i=1}^n (\frac{w_{i}^{(n)}}{n} -\frac{1}{n}  ) X_i}{S_n \sqrt{\sum_{i=1}^n (\frac{w_{i}^{(n)}}{n} -\frac{1}{n}  )^2} }
\end{equation}
to the same context will also be explored (cf. Section \ref{Bootstrap and Linear Processes}).

\par
Just like in \cite{Csorgo and Nasari}, in this paper the  method  of conditioning on the random  weights $w_{i}^{(n)}$'s   is used for constructing randomized pivots. Viewing the randomized sums of linear processes  as weighted sums of the original data, here we derive the asymptotic normality of properly normalized/Studentized   randomized  sums of short and long memory linear processes.   As will be seen,  our conditional CLTs  also   imply unconditional CLTs in terms of  the joint distribution of the observables and the random weights.

%\par
%Another interesting feature of  this paper is the use of stochastically weighing $(X_i -\mu)$ in $G_{n}$, as defined in  (\ref{eq 6}), to establish  CLTs  for certain linear processes for which the use of their variance as a normalizing scalar leads to convergence to  a degenerate limit (cf. Remark \ref{overdifferencing}). %The same is also  true for the Studentized versions of the latter mentioned statistics.
\par
The material in this paper is organized as follows. In Section \ref{Main Results} the randomized pivots are introduced and  conditional and  unconditional CLTs are presented for them.  In Section \ref{Confidence intervals}, asymptotic confidence intervals of size $1-\alpha$, $0< \alpha<1$, are constructed for the population mean $\mu=E X_1$. Also in Section \ref{Confidence intervals},  confidence bounds  are constructed  for some functionals of linear processes. The results in this section are directly applicable to constructing   confidence bounds for some  functionals of long memory  linear processes whose limiting distribution may not  necessarily be normal.  Section \ref{Simulation Results} is devoted to presenting  our simulations and numerical studies. In Section \ref{Bootstrap and Linear Processes}, we study the problem of bootstrapping linear processes and provide a comparison between our results and those obtained by using the bootstrap.  The proofs are given in Section \ref{Proofs}.

\section{CLT for randomized  pivots of the population mean }\label{Main Results}
Throughout  this section,  we let $\{X_i; i\geq 1\}$ be a linear process that, for each $i \geq 1$,
 is   defined by
\begin{equation}\label{eq 1}
X_i=\mu+\sum_{k=0}^{\infty} a_{k} \zeta_{i-k}=\mu+\sum_{k=-\infty}^{i} a_{i-k}\zeta_{k},
\end{equation}
where $\mu $ is a real number, $\{a_k; k \in \mathbb{Z}  \}$ is a sequence of real numbers such that $\sum_{k=0}^{\infty} a^{2}_{k}<+\infty$ and $\{\zeta_k; k \in \mathbb{Z} \}$  are i.i.d. white noise innovations  with $E \zeta_k=0$ and $0<\sigma_{\zeta}^2:=Var(\zeta_k)<+\infty$. Consequently,   we have
$E X_i=\mu$. Moreover,  we assume throughout  that the    $X_i$s, are non-degenerate and  have a  finite variance    $
\gamma_{0}:= E X_{i}^2-\mu^2:=\sigma_{\zeta}^{2} \sum_{k=0}^{\infty} a^{2}_{k}-\mu^2$, $i \geq 1$. We note in passing that for some of the results in this paper  the existence and finiteness of some  higher moments of the data  will also  be assumed (cf. Theorems \ref{CLT G^*^stu} and \ref{CLT T^*^stu}).
\par
For throughout  use, we let
\begin{equation}\label{eq 1'}
\gamma_h:= Cov(X_s,X_{s+h})=E(X_{s}-\mu)(X_{s+h}-\mu), \ h\geq 0, \ s\geq 1,
\end{equation}
be the autocovariance function of the stationary linear process $\{X_i, i\geq 1 \}$ as in (\ref{eq 1}).
 %with the spectral density of the form
%\begin{equation}\label{spectral density}
%f(\lambda)=b_0 |\lambda|^{-2d}+o(|\lambda|^{-2d}), \ as \ \lambda\to 0,
%\end{equation}
%where   $b_0>0$ and $d$, $0\leq d <1/2$, is the memory parameter.
Moreover, based on the stationary sample $X_1, \ldots,X_n$,   $n\geq 1$, on the linear process $\{X_i, i \geq 1 \}$,  for throughout use we define
\begin{eqnarray}
&&\bar{X}_{n}:=\sum_{i=1}^{n}X_i\big/n, \nonumber\\
&&\bar{\gamma}_i:=\sum_{j=1}^{n-i}(X_j -\bar{X}_{n})(X_{j+i}-\bar{X}_n)\big/ n, \ 0\leq i\leq n-1, \label{eq 2}
\end{eqnarray}
respectively  the sample mean and  sample autocovariance.

% \par
%We note  that in this paper we consider linear processes, as in (\ref{eq 1}), with $0\leq d<1/2 $, where the case  $d=0$ represents a \emph{short memory} process, and in the case of $0<d<1/2$,  the  process   exhibits  a pattern of  \emph{long memory}. Furthermore, we  note that the assumption that the linear process in hand has a spectral density  of the form (\ref{spectral density}) is needed only for part (B) of Theorem  \ref{CLT G^*^stu} and it is not required for the conclusions of Theorem  \ref{CLT G^*} to hold true.

\par
We now define the following two randomly weighted versions of the partial sum  $\sum_{i=1}^{n} (X_{i}-\mu)$, where $\{X_i, i \geq 1\}$ is a sequence of linear processes defined in (\ref{eq 1}),
\begin{eqnarray}
&&\sum_{i=1}^{n} (\frac{w_{i}^{(n)}}{n}-\frac{1}{n})(X_i-\mu),  \label{eq 3} \\
&&\sum_{i=1}^{n} \big| \frac{w_{i}^{(n)}}{n}-\frac{1}{n} \big|(X_i -\mu),\label{eq 4}
\end{eqnarray}
where  the random  weights in the triangular array  $\{w^{(n)}_{1},\ldots,w^{(n)}_{n}\}_{n=1}^{\infty}$ have a multinomial distribution of size
$n= \sum_{i=1}^n  w^{(n)}_{i}$ with respective probabilities  $1/n$, i.e.,
\begin{equation*}
(w^{(n)}_{1},\ldots,w^{(n)}_{n})\  \substack{d\\=}\ \ multinomial(n;\frac{1}{n},\ldots,\frac{1}{n}),
\end{equation*}
are independent from the stationary sample $\{ X_1,\ldots,X_n\}$, $n\geq 1$, on the process $\{X_{i}, i \geq 1 \}$ as in (\ref{eq 1}).
\par
The just introduced  randomized sums in  (\ref{eq 3}) and  (\ref{eq 4}), which are randomized  versions of $\sum_{i=1}^n (X_i-\mu)$, can be computed via generating  a realization of the multinomial random weights $(w_{1}^{(n)},\ldots,w_{n}^{(n)})$ with $\sum_{i=1}^n w_{i}^{(n)}=n$ with associated probability vector $(1/n,\ldots,1/n)$. In this context, one  way of generating the random weights $w_i^{(n)}$'s is to
re-sample  from the set of indices $\{1,\ldots,n \}$ of the stationary sample  $X_1,\ldots,X_n$ in hand, $n\geq 1$,  with replacement $n$ times with respective probabilities $1/n$ so that,  for each $1\leq i \leq n$, $w^{(n)}_i$ is the  number of times the index $i$ of $X_i$ is chosen in this re-sampling process.

\begin{remark}\label{Remark 0}
In view of $ \sum_{i=1}^n  w^{(n)}_{i}=n$, one can readily see that for the randomized sum defined in (\ref{eq 3}), we have
\begin{eqnarray*}
\sum_{i=1}^{n} (\frac{w_{i}^{(n)}}{n}-\frac{1}{n})(X_i-\mu)&=& \sum_{i=1}^{n} (\frac{w_{i}^{(n)}}{n}-\frac{1}{n}) X_i\\
&=:& \bar{X}_{n}^{*} - \bar{X}_n,
\end{eqnarray*}
i.e., $\bar{X}_{n}^{*} - \bar{X}_n$, which coincides with    the randomized sum  (\ref{eq 3}), forgets about what the value of the population mean  $\mu=E X_1$ might  be. On the other hand,    the randomization used in  the sum (\ref{eq 4}) preserves $\mu=E X_1$.
\end{remark}

\par
In addition to  preserving $\mu$, the  randomized sum (\ref{eq 4}) tends  to preserve the covariance structure of the data as well, a property that the sum (\ref{eq 3}) fails to maintain (cf. Remark \ref{T^* no good long memory}).

\par
Properly normalized,     (\ref{eq 4})  provides a  natural direct  pivot for the population mean $\mu=E X_1$ (cf. the definition (\ref{eq 6})).

\par
For throughout  use we introduce the following notations.
\\
\textbf{Notations}. Let $(\Omega_X,\mathfrak{F}_X,P_X)$ denote  the probability space of the random variables  $X,X_1,\ldots$, and $(\Omega_w,\mathfrak{F}_w,P_w)$ be  the probability space on which $\big(w^{(1)}_1,(w^{(2)}_1,w^{(2)}_{2}),\ldots,(w^{(n)}_1,\ldots,w^{(n)}_{n}),\ldots \big)$ are defined. In view of the independence of these two sets of random variables, jointly they live on the direct product probability space $(\Omega_X \times \Omega_w, \mathfrak{F}_X \otimes \mathfrak{F}_w, P_{X,w}=P_X \ .\ P_w )$. For each $n\geq 1$, we also let  $P_{.|w}(.)$   stand for the conditional probability   given $\mathfrak{F}^{(n)}_w:=\sigma(w^{(n)}_1,\ldots,w^{(n)}_{n})$  with corresponding conditional expected value and variance,   $E_{.|w}(.)$ and $Var_{.|w}(.)$, respectively.
\par
In a similar fashion to the randomized pivots in the i.i.d. case as  in   (\ref{added Cs2}), in this context  we define the randomized $t$-type statistics, based on   the sample $X_1,\ldots,X_n$, $n\geq 1$, of linear processes defined by (\ref{eq 1}),  as follows
\begin{eqnarray}
G_{n}&:=& \frac{\sum_{i=1}^{n} \big|\frac{w_{i}^{(n)}}{n}-\frac{1}{n}\big| (X_i-\mu)}{\sqrt{\gamma_0 \sum_{j=1}^{n} (\frac{w_{i}^{(n)}}{n}-\frac{1}{n})^2+ 2\sum_{h=1}^{n-1} \gamma_{h} \sum_{j=1}^{n-h} \big|\frac{w_{j}^{(n)}}{n}-\frac{1}{n}\big|  \    \big|\frac{w_{j+h}^{(n)}}{n}-\frac{1}{n}\big| } },\nonumber\\
&&\label{eq 6}
\end{eqnarray}
where $\gamma_h$, $0\leq h \leq n-1$, are defined in (\ref{eq 1'}).
\par
We note that the denominator  of  $G_{n}$ is
\begin{equation*}
\Big(Var_{X|w}\big( \sum_{i=1}^{n} \big|\frac{w_{i}^{(n)}}{n}-\frac{1}{n}\big| (X_i-\mu) \big)\Big)^{1/2}.
\end{equation*}

\par
Noting that the normalizing sequence  in   $G_{n}$  depends on the  parameters $\gamma_h$, $0\leq h \leq n-1$, we now define  Studentized  versions  of it.
\par
The following Studentized statistic,  in (\ref{eq 8}), is defined for  all short memory  linear processes as in (\ref{eq 1}) as well as, for   long memory linear processes as defined in (\ref{eq 1}) with    $a_k \sim c k^{d-1}$, for some $c>0$, as $k \rightarrow +\infty$, where $0< d <1/2$. We refer to $d$ as the memory parameter. In order to    unify our notation in (\ref{eq 8}) to include both short and long memory linear processes, when the linear process in hand is of short memory,  then  we define the memory parameter   $d$ to be zero. Thus, the Studentized version of $G_n$ as in (\ref{eq 6}) is defined as follows

\begin{equation}\label{eq 8}
G_{n}^{{stu}}(d):= \frac{\sum_{i=1}^{n} \big|\frac{w_{i}^{(n)}}{n}-\frac{1}{n}\big| (X_i-\mu)}{\sqrt{(\frac{q}{n})^{-2d}\bar{\gamma}_0 \sum_{j=1}^{n} (\frac{w_{i}^{(n)}}{n}-\frac{1}{n})^2+ 2 \sum_{h=1}^{q} \bar{\gamma}_{h} \sum_{j=1}^{q-h} \big|\frac{w_{j}^{(n)}}{n^{1-2d}}-\frac{1}{n^{1-2d}}\big|    \big|\frac{w_{j+h}^{(n)}}{q^{1+2d}}-\frac{1}{q^{1+2d}}\big| } },
\end{equation}
where $q \to +\infty$  in such a way that as $n\to +\infty$, $q=O(n^{1/2})$,  $\bar{\gamma}_{h}$ is  defined    in (\ref{eq 2}) and $0\leq d<1/2$.

\par
We note in passing that  in the case of long memory linear processes, when an estimator $\hat{d}$ is used to replace the memory parameter $d$ (cf. Section \ref{Simulation Results}),    then the  notation $G_{n}^{{stu}}(\hat{d})$ stands for the version of $G_{n}^{{stu}}(d)$, as in (\ref{eq 8}),  in which $d$ is replaced by $\hat{d}$. Also, in the case of having a short memory linear process, i.e, when $d=0$, the notation  $G_{n}^{{stu}}(0)$ stands for the version of $G_{n}^{{stu}}(d)$ in which $d$ is replaced by 0.
\par
It can also  be readily seen that, in the case of having a long memory linear process,  after estimating $d$ by a proper estimator  $\hat{d}$, then $G_{n}^{{stu}}(\hat{d})$, apart from $\mu$ that is to be estimated,   is  computable based on the data $X_1,\ldots,X_n$ and the generated  multinomial weights $(w^{(n)}_{1},\ldots,w^{(n)}_n)$. The same is also true when dealing with short memory linear processes, i.e., when $d=0$. In other words, in the case of short memory linear processes, apart from the population mean $\mu$, which is to be estimated, the other elements of the pivot   $G_{n}^{{stu}}(0)$ are computable based on the data and the generated  multinomial  weights.

\par
The following two theorems, namely Theorems \ref{CLT G^*} and \ref{CLT G^*^stu},  establish conditional (given the weights) and unconditional  CLTs  for $G_{n}$ and $G_{n}^{{stu}}(d)$, respectively.
 These theorems  are valid for classes of  both  short  and long memory data.
 \par
We note that throughout this paper $\Phi(.)$ stands  for the standard normal distribution function.

\begin{thm}\label{CLT G^*}
Suppose that $\{X_i, i \geq 1 \}$ is a stationary linear process as defined  in (\ref{eq 1}) with $\sum_{k=0}^{\infty} a^{2}_{k}<+\infty$.
\\
(A) If for each $h\geq 1$,  $ \gamma_{h}\geq 0$, then,  as $n \to +\infty$, we have for all $t\in \mathbb{R}$,
\begin{equation}\label{eq 13}
P_{X|w} (G_{n}\leq t) \longrightarrow \Phi(t)\ in \ probability-P_w
\end{equation}
and, consequently,
\begin{equation}\label{eq 14}
P_{X,w} (G_{n}\leq t) \longrightarrow \Phi(t), \ t\in \mathbb{R}.
\end{equation}
(B) If  $\sum_{k=0}^{\infty} |a_k|<+\infty$ and
\begin{equation}\label{eq 15}
 \gamma_0 +8 e^{-2}    \sum_{h= 1}^{+\infty}  \gamma_h \neq 0,
\end{equation}
then, for all $t\in \mathbb{R}$, as $n \to +\infty$,  the conditional CLT (\ref{eq 13}) and, consequently, also the unconditional CLT (\ref{eq 14})  hold true.
\end{thm}

\par
We  note that part (A) of Theorem \ref{CLT G^*} holds true for all kinds of stationary short and long memory linear processes with $\gamma_h\geq 0$, $h\geq 1$. For short memory processes with possibly some negative $\gamma_h$s,  we continue to  have  (\ref{eq 13}) and its consequence (\ref{eq 14}), as long as (\ref{eq 15}) holds true.

\begin{remark}\label{overdifferencing}
 Theorem  \ref{CLT G^*}  allows having  CLTs, conditionally on the weights,  or in terms of the joint distribution of the data and the random weights,  for randomized   versions of  partial sums of   linear processes for which there are no CLTs  with the standard deviation of the partial sum in hand in its normalizing sequence. Examples  of such processes, which are usually the  results of \emph{overdifferencing}, are of the form  $X_t=Y_t -Y_{t-1}$, where the $Y_t$ are white noise,  like, e.g.,   the well known non-invertible moving average MA(1) processes.  Randomizing these processes   results in  randomly weighted  partial sums of the original data whose  variance,  unlike the variance of the original partial sums,  go  to infinity  as the sample size $n\to +\infty$. This phenomenon can be seen to be the result of incorporating the random weights, for then  the sum  $\sum_{i=1}^n \big|  \frac{w_{i}^{(n)}}{n}-\frac{1}{n} \big| X_i$ no longer forms a telescoping  series as the original non-randomized sum $\sum_{t=1}^n X_t= Y_n -Y_0$.
%We  illustrate  this phenomenon  using   simulations  in Section \ref{Simulation Results}.
\end{remark}

\par
The following  result, which is a companion  of Theorem  \ref{CLT G^*}, establishes the asymptotic normality for the Studentized statistics  $G_{n}^{{stu}}(d)$.

\begin{thm}\label{CLT G^*^stu}
(A) Assume that    the stationary  linear process $\{X_i, i \geq 1 \}$, as defined in (\ref{eq 1}), is of  short memory, i.e., $\sum_{k=0}^{\infty} |a_{k}|<+\infty$,  and   $E \zeta_{1}^{4}<+\infty$. Also, assume that     $\gamma_{h}\geq 0$, for all $h\geq 1$.
 Then, as $n,q\to +\infty$ such that $q=O(n^{1/2})$,   we have for all $t \in \mathbb{R}$,
\begin{equation}\nonumber\label{eq 18}
P_{X|w} (G_{n}^{{stu}}(0)\leq t) \longrightarrow \Phi(t) \ in \ probability-P_{w}
\end{equation}
and, consequently,
\begin{equation}\nonumber\label{eq 19}
P_{X,w} (G_{n}^{{stu}}(0)\leq t) \longrightarrow \Phi(t),\ t \in \mathbb{R}.
\end{equation}
(B) Let  the linear process $\{X_i, i \geq 1 \}$, as defined in (\ref{eq 1}), with $\sum_{k=0}^{\infty} a^{2}_{k}<+\infty$ and $\gamma_{h}\geq 0$, for all $h\geq 1$,  be of long memory  such that $E \zeta_{1}^{4}<+\infty$ and, as $k \rightarrow +\infty$,  $a_k \sim c k^{d-1}$, for some $c>0$, where $0< d <1/2$. Then, as $n,q\to +\infty$ such that $q=O(n^{1/2})$, for all $t \in \mathbb{R}$, we have

\begin{eqnarray*}
&&P_{X|w} (G_{n}^{{stu}}(d)\leq t) \longrightarrow \Phi(t) \ in \ probability-P_{w}\\
&&P_{X|w} (G_{n}^{{stu}}(\hat{d})\leq t) \longrightarrow \Phi(t) \ in \ probability-P_{w},
\end{eqnarray*}
and, consequently,

\begin{eqnarray*}
&&P_{X,w} (G_{n}^{{stu}}(d)\leq t) \longrightarrow \Phi(t),\ t \in \mathds{R} \\
&&P_{X,w} (G_{n}^{{stu}}(\hat{d})\leq t) \longrightarrow \Phi(t), \ t \in \mathds{R},
\end{eqnarray*}
where $\hat{d}$ is an estimator of the memory parameter $d$ such that $\hat{d}-d=o_{P_{X}}(1/\log n)$.
\end{thm}

\section{Randomized  confidence intervals for the population mean $\mu$ }\label{Confidence intervals}
In this section, we  use  $G^{stu}_n(\hat{d})$  as a natural randomized pivot for the population mean $\mu$ in a nonparametric way. Based on it, we now spell out   asymptotic randomized $1-\alpha$ size  confidence intervals for the  population mean  $\mu$. In what  follows $z_{1-\alpha}$ stands for the solution to
$\Phi( z_{1-\alpha})=1-\alpha$. For the ease of  notation, we first introduce  the following setup:
\begin{eqnarray*}\nonumber
D_{n,q,0}&:=&  \bar{\gamma}_{0} \sum_{j=1}^{n} (\frac{w_{i}^{(n)}}{n}-\frac{1}{n})^2+ 2 \sum_{h=1}^{q} \bar{\gamma}_{h} \sum_{j=1}^{q-h} \big|\frac{w_{j}^{(n)}}{n}-\frac{1}{n}\big|     \big|\frac{w_{j+h}^{(n)}}{q}-\frac{1}{q}\big|, \\
D_{n,q,\hat{d}}&:=& (\frac{q}{n})^{-2\hat{d}}
\bar{\gamma}_{0} \sum_{j=1}^{n} (\frac{w_{i}^{(n)}}{n}-\frac{1}{n})^2+ 2 \sum_{h=1}^{q} \bar{\gamma}_{h} \sum_{j=1}^{q-h} \big|\frac{w_{j}^{(n)}}{n^{1-2\hat{d}}}-\frac{1}{n^{1-2\hat{d}}}\big|    \big|\frac{w_{j+h}^{(n)}}{q^{1+2\hat{d}}}-\frac{1}{q^{1+2\hat{d}}}\big| .
\end{eqnarray*}

It is important to note that the randomized confidence intervals (one or two-sided) which we are about to present,   henceforth,  are valid in terms of the conditional distribution $P_{X|w}$,  as well as in terms of the joint distribution $P_{X,w}$.

\par
When the linear process in hand possesses the property of   short memory,  if it satisfies the conditions of part (A) of Theorem \ref{CLT G^*^stu}, then  the  asymptotic two-sided  $1-\alpha$ size randomized  confidence interval for the population mean $\mu=E_{X} X_1$ has the following form.
\begin{equation}\label{C. I. for mu2}
\frac{\sum_{i=1}^n \big|\frac{w_{i}^{(n)}}{n}-\frac{1}{n} \big| X_i-  z_{1-\alpha/2} D^{1/2}_{n,q,0} }{\sum_{j=1}^n \big|\frac{w_{j}^{(n)}}{n}-\frac{1}{n} \big|}\leq  \mu  \leq \frac{\sum_{i=1}^n \big|\frac{w_{i}^{(n)}}{n}-\frac{1}{n} \big| X_i +  z_{1-\alpha/2} D^{1/2}_{n,q,0} }{\sum_{j=1}^n \big|\frac{w_{j}^{(n)}}{n}-\frac{1}{n} \big|}
\end{equation}

\par
An asymptotic $1-\alpha$ size randomized  two-sided confidence interval for the population mean  $\mu=E_{X} X_1$ of a long range dependent linear process,   as defined in (\ref{eq 1}), when it satisfies  the   conditions in part  (B) of Theorem \ref{CLT G^*^stu} is constructed as follows.

\begin{equation}\label{C. I. for mu1}
\frac{ \sum_{i=1}^n \big|\frac{w_{i}^{(n)}}{n}-\frac{1}{n} \big| X_i  -  z_{1-\alpha/2} D^{1/2}_{n,q,\hat{d}} }{\sum_{j=1}^n \big|\frac{w_{j}^{(n)}}{n}-\frac{1}{n} \big|}\leq \mu  \leq \frac{\sum_{i=1}^n \big|\frac{w_{i}^{(n)}}{n}-\frac{1}{n} \big| X_i +  z_{1-\alpha/2} D^{1/2}_{n,q,\hat{d}} }{\sum_{j=1}^n \big|\frac{w_{j}^{(n)}}{n}-\frac{1}{n} \big|}
\end{equation}

\subsection{Confidence bounds for the mean of some functionals of long memory linear processes}\label{Functionals}
\par
The following result, namely Corollary \ref{corollary 1},  is a consequence of   Theorem \ref{CLT G^*^stu} and  Jensen's inequality. Corollary \ref{corollary 1}     gives randomized confidence bounds  for $\mu_{\mathcal{G}}:=E_{X}\mathcal{G}(X_i)$, for some measurable functions $\mathcal{G}$,  i.e., for the mean of certain subordinated  functions of  the long memory   linear process  in hand,   and it reads as follows.

\begin{corollary}\label{corollary 1}
Let $\{X_i, i \geq 1\}$ be so that it satisfies the conditions in  (B) of Theorem \ref{CLT G^*^stu}.   Assume that $\mu_\mathcal{G}=E_{X} | \mathcal{G}(X_i) | <+\infty$.
As $n,q\to +\infty$ in such a way that $q=O(n^{1/2})$, we have
\\
\\
(A) If $\mathcal{G}$ is increasing and convex, then, an asymptotic  $1-\alpha$ size lower confidence bound for $\mu_{G}$ is
\begin{equation}\nonumber
\mu_{\mathcal{G}}\geq \mathcal{G}(\frac{\sum_{i=1}^n \big|\frac{w^{(n)}_{i}}{n}-\frac{1}{n} \big|   X_i- z_{1-\alpha} D^{1/2}_{n,q,\hat{d}}} {\sum_{i=1}^n \big|\frac{w^{(n)}_{i}}{n}-\frac{1}{n} \big|} )
\end{equation}
(B) If $\mathcal{G}$ is decreasing  and convex, then, a randomized  asymptotic  $1-\alpha$ size lower confidence bound for $\mu_{\mathcal{G}}$ is
\begin{equation}\nonumber
\mu_{\mathcal{G}}\geq \mathcal{G}(\frac{\sum_{i=1}^n \big|\frac{w^{(n)}_{i}}{n}-\frac{1}{n} \big|   X_i+ z_{1-\alpha} D^{1/2}_{n,q,\hat{d}}} {\sum_{i=1}^n \big|\frac{w^{(n)}_{i}}{n}-\frac{1}{n} \big|} ).
\end{equation}

\end{corollary}

\begin{remark}
 Corollary  \ref{corollary 1} remains valid for functionals of  short memory linear processes with $D_{n,q,0}$ replacing  $D_{n,q,\hat{d}}$.   It is also important to note that  the conclusions of (A) and  (B) of Corollary  \ref{corollary 1} hold true without making any assumptions   about the variance of the subordinated function $\mathcal{G}$. In other words, Corollary \ref{corollary 1} is valid even when $Var({\mathcal{G}})$  is  not finite.
\end{remark}

\begin{remark}
In reference to studying the mean of functions of stationary long memory Gaussian linear processes, Corollary \ref{corollary 1} helps  avoiding   dealing with the  sampling distributions of  processes of the  functions of  stationary long memory Gaussian processes which are known to be relatively complicated, specially when they exhibit non-normal asymptotic distributions (cf. Taqqu \cite{Taqqu} and Dobrushin  and  Major \cite{Dobrushin  and  Major}, for example). We note that any long memory Gaussian process $\{\eta_i;\ i\geq 1 \}$, i.e., $Cov(\eta_1,\eta_{1+k})=k^{-\alpha}\ L(k)$, where $0<\alpha <1$ and $L(.)$ is a slowly varying function at infinity, can be viewed as a long memory  linear process (cf. Cs\'{a}ki  \emph{et al}. \cite{Csaki et al}, for example) that satisfies the conditions of part (B) of Theorem \ref{CLT G^*^stu}. Therefore,  Corollary \ref{corollary 1} is directly applicable to constructing randomized confidence bounds for  means of subordinated functions of long memory Gaussian processes $\{\eta_i; \ i\geq 1\}$ without    making assumptions concerning their   variance,   or   referring to their  Hermit expansions.
\end{remark}

\section{Simulation results}\label{Simulation Results}
In this section we examine numerically   the performance  of   $G^{{stu}}_n (d)$ and $G^{{stu}}_n (\hat{d})$, in view of the CLTs in  Theorem  \ref{CLT G^*^stu}, versus those  of their classical counterparts $T^{stu}_{n} (d)$ and $T^{stu}_{n} (\hat{d})$, defined as
\begin{equation}\label{Classic T with d}
T^{stu}_{n} (d):= \frac{n^{1/2 -d} (  \bar{X}_n -\mu) }{\sqrt{ q^{-2d} \bar{\gamma}_0  + 2 q^{-2d} \sum_{h=1}^{q} \bar{\gamma}_{h}  (1-h/q)  } },
\end{equation}

\begin{equation}\label{Classic T with dhat}
T^{stu}_{n} (\hat{d}):= \frac{n^{1/2 -\hat{d}} (  \bar{X}_n -\mu) }{\sqrt{ q^{-2\hat{d}} \bar{\gamma}_0  + 2 q^{-2\hat{d}} \sum_{h=1}^{q} \bar{\gamma}_{h}  (1-h/q)  } }.
\end{equation}
When the linear process in hand is of short memory, then $T^{stu}_{n} (0)$ stands for a version of  $T^{stu}_{n} (d)$ in which $d$ is replaced by 0.
Under the conditions of  our Theorem \ref{CLT G^*^stu}, from Theorem 3.1 of Graitis \emph{et al}. \cite{Giraitis et al.} and Theorem 2.2 of Abadir \emph{et al}. \cite{Abadir  Distaso Giraitis},  we conclude that the limiting distribution of  $T^{stu}_{n} (0)$, $T^{stu}_{n} (d)$ and $T^{stu}_{n} (\hat{d})$ is standard normal. Hence, $T^{stu}_{n} (0)$, $T^{stu}_{n} (d)$ and $T^{stu}_{n} (\hat{d})$ converge to the same limiting distribution as that of $G^{{stu}}_n (0)$, $G^{{stu}}_n (d)$ and $G^{{stu}}_n (\hat{d})$, under the conditions  of Theorem \ref{CLT G^*^stu}.

\par
In Tables 1-6, we provide  motivating simulation results in preparation for the upcoming     in depth numerical studies in Tables 7-12. In the following Tables 1-6 we use packages ``arima.sim'' and ``fracdiff.sim'' in R    to  generate observations from short and  long memory \emph{Gaussian} processes, respectively.  Tables 1-6 present empirical probabilities of  coverage   with the normal  cutoff points $\pm 1.96$, i.e.,  the nominal probability coverage is $0.95$. The results are based on  1000 replications of the data and the multinomial weights $(w_{1}^{(n)},\ldots,w_{n}^{(n)})$.   The choice of $q$ is made based on relation (2.14) of  Abadir \emph{et al}. \cite{Abadir  Distaso Giraitis} in each case. More precisely,  we let  $q$ be $ceiling(n^{1/3})$ for the examined  short memory linear processes, and for long memory linear processes with $0< d <0.25$,   we let $q$ be $ceiling(n^{1/(3+4d)})$,  and for $0.25<d<0.5$, we let $q$ be    $ceiling(n^{1/2 -d})$ and, when the data are long memory with parameter $d$, then $\hat{d}$ stands for the MLE approximation    of $d$, with the   Haslett and Raftery  \cite{Haslett and Raftery} method used to approximate the likelihood. This estimator of $d$ is provided in the R package ``fracdiff'' and it is  used in our simulation studies in Tables 4, 6, 10 and 12. We note that there are other commonly used methods of estimating the  memory parameter $d$,  such as the  Whittle estimator (cf. K\"{u}nch  \cite{K\"{u}nch} and Robinson \cite{Robinson}), which is available in the R package ``longmemo''  using  the   Beran  \cite{Beran}  algorithm.
For more on estimators for    the  memory parameter $d$ and their asymptotic behavior, we refer to, for example, Robinson  \cite{Robinson 1997} and   Mulines and Soulier  \cite{Moulines and Soulier} and references therein.

\newpage
\begin{table}
\caption{MA(1): \ $ X_{t}=W_t-0.5\ W_{t-1}$  }
\begin{center}
%Model: $ X_{t}$ is $ARIMA(0,0,0.5) $\\
\begin{tabular}{ c|c|c|c   }
\hline \hline
  Distribution &  $n$ & $\textmd{Coveage\ prob.\ of\ }  G_{n}^{{stu}}(0)$ & $\textmd{Covergae prob. of\ } T_{n}^{stu}(0)$  \\
\hline\hline
 \multirow{2}{*}{$ W_t \ {\substack{d\\=}}\ \ \textrm{N(0,1)}$} %&15  & 0.926 & 0.915 \\
                                  &20  & 0.936 & 0.933   \\
                                  &30  & 0.949 & 0.943 \\
%                         \hline \hline
% \multirow{2}{*}{$W_t \ {\substack{d\\=}}\ \ \textrm{lognormal(0,1)}$} &20&  0.160 & 0.004 \\
 %                                   &30&  0.210 & 0.156 \\
\hline
\end{tabular}

\end{center}
 \end{table}

\begin{table}
\caption{AR(1):\ $X_t=0.5 \ X_{t-1}+W_t$ }
\begin{center}
%Model: $X_t$ is $ ARIMA(0.5,0,0)$\\
\begin{tabular}{ c|c|c|c   }
\hline \hline
  Distribution &  $n$ &  $\textmd{Coverage \ prob.\ of}\ G_{n}^{{stu}}(0)$ & $\textmd{Coverage \ prob. \ of} \ T_{n}^{stu}(0)$  \\
\hline\hline
 \multirow{2}{*}{$W_t\ {\substack{d\\=}}\ N(0,1)$} &20  & 0.935 & 0.927 \\
                                  &30 & 0.947 & 0.941 \\
                                  %&50 & 0.957 & 0.965\\
%                         \hline \hline
% \multirow{2}{*}{$W_t\ {\substack{d\\=}}\ lognormal(0,1)$} &70& 0.302 & 0.130 \\
 %                                   &80&  0.536 & 0.384 \\

\hline
\end{tabular}

\end{center}
\end{table}

\begin{table}[H]\label{L1}
\caption{Long Memory with $d=0.2$:\ $X_{t}=(1-B)^{0.2}\ W_t$ }
\begin{center}
%Model: $ X_{t}$ is long memory with $d=0.2$\\
\begin{tabular}{ c|c|c|c   }
\hline \hline
  Distribution &  $n$ & $\textmd{Coverage \ prob.\ of}\ G_{n}^{{stu}}(0.2)$ &  $\textmd{Coverage \ prob.\ of}\ T_{n}^{stu} (0.2)$\\
\hline\hline
 \multirow{2}{*}{$W_t\ {\substack{d\\=}}\ N(0,1)$} &30  & 0.933 &0.882 \\
                                  &50  & 0.956 & 0.921  \\
                                 % &300  & 0.952 0.899    \\

%                         \hline \hline
% \multirow{2}{*}{$W_t \ {\substack{d\\=}} \ lognormal(0,1)$} &150  &  0.588 &  0.004  \\
%                                    &250  & 0.672 &  0.022 \\
\hline
\end{tabular}

\end{center}
\end{table}

\begin{table}[H]\label{Long2}
\caption{Long Memory with $d=0.2$:\ $X_{t}=(1-B)^{0.2} \ W_t$; estimator $\hat{d}$ used}
\begin{center}
\begin{tabular}{ c|c|c|c }
\hline \hline
  Distribution &  $n$ & $\textmd{Coverage \ prob.\ of}\ G_{n}^{{stu}}(\hat{d})$ &  $\textmd{Coverage \ prob.\ of}\ T_{n}^{stu} (\hat{d})$\\
\hline\hline
 \multirow{2}{*}{$W_t\ {\substack{d\\=}}\ N(0,1)$} &200  & 0.939 &0.887  \\
                                                   &300  & 0.945& 0.893   \\
                                                   %&400  & 0.950 & 0.904 \\
%                         \hline \hline
 %\multirow{2}{*}{$W_t\ {\substack{d\\=}}\ lognormal(0,1)$} &400  & 0.584 & 0.002  \\
                                   % &500  & 0.656 &  0.002  \\
\hline
\end{tabular}

\end{center}
\end{table}

\begin{table}[H]\label{Long3}
\caption{Long Memory with $d=0.4$:\ $X_{t}=(1-B)^{0.4}\ W_t$}
\begin{center}
%Model: $ X_{t}$ is long memory with $d=0.4$\\
\begin{tabular}{ c|c|c|c }
\hline \hline
  Distribution &  $n$ & $\textmd{Coverage \ prob.\ of}\ G_{n}^{{stu}}(0.4)$ &  $\textmd{Coverage \ prob.\ of}\ T_{n}^{stu} (0.4)$\\
\hline\hline

\multirow{2}{*}{$W_t\ {\substack{d\\=}}\ N(0,1)$}  &300  & 0.938& 0.891 \\
                                                   &400  & 0.946& 0.903 \\
                                                   %&200  & 0.942 & 0.893 \\
                                 %&300  & 0.955 &  0.908 \\
 %                        \hline \hline
% \multirow{2}{*}{$W_t\ {\substack{d\\=}}\ lognormal(0,1)$} &300  & 0.148 & 0.000 \\
%                                    &400  & 0.328 & 0.000  \\
\hline
\end{tabular}

\end{center}
\end{table}

\begin{table}[H]\label{G4}
\caption{Long Memory with $d=0.4$:\ $X_{t}=(1-B)^{0.4}\ W_t$; estimator $\hat{d}$ used}
\begin{center}
%Model: $ X_{t}$ is long memory with $d=0.4$\\
\begin{tabular}{ c|c|c|c }
\hline \hline
  Distribution &  $n$ & $\textmd{Coverage \ prob.\ of}\ G_{n}^{{stu}}(\hat{d})$ &  $\textmd{Coverage \ prob.\ of}\ T_{n}^{stu} (\hat{d})$\\
\hline\hline
 %\multirow{2}{*}{$W_t\  {\substack{d\\=}}\ N(0,1)$} &800  & 0.150 &0.000 \\
%                                  &1000  & 0.328 &0.000 \\
\multirow{2}{*}{$W_t\  {\substack{d\\=}}\ N(0,1)$} &500   & 0.882 & 0.842 \\
                                                   &1000  &  0.926 & 0.876  \\

\hline
\end{tabular}

\end{center}
\end{table}

\par
We  now  present a more in depth simulation  study   for   non-Gaussian linear processes. Here, once again,   we use the packages ``arima.sim'' and ``fracdiff.sim'' in R    to  generate observations from short and  long memory linear processes, respectively. In our numerical studies below,  we use the standardized   Lognormal(0,1) distribution, i.e., Lognormal with mean 0 and variance 1,  to generate observations from  short or long memory linear  processes.
The choice of Lognormal(0,1) in our studies is due to the fact it is a heavily  skewed distribution.

\par
The following Tables 7-12 are presented to illustrate the significantly better  performance  of  $G_{n}^{stu}(d)$, $G_{n}^{stu}(\hat{d})$ and $G_{n}^{stu}(0)$  over their respective classical counterparts $T_{n}^{stu}(d)$, $T_{n}^{stu}(\hat{d})$ and $T_{n}^{stu}(0)$,   in view of   our Theorem \ref{CLT G^*^stu}.
  In  Tables 7 and 8  we present  numerical comparisons between the performance of  $G^{{stu}}_{n}(0)$ to that of $T^{{stu}}_{n}(0)$, both as pivots for the population mean,  for some moving average and autoregressive processes.
 The numerical performance of $G^{{stu}}_{n}(d)$ to that of the  classical $T^{{stu}}_{n}(d)$ for some long memory linear processes  are presented in Tables 9 and 11. Tables 10 and 12 are specified to comparing $G^{{stu}}_{n}(\hat{d})$ to $T^{{stu}}_{n}(\hat{d})$.
\par
In Table 7, for the therein underlined MA(1) process,  we generate 500 empirical coverage probabilities of the event that
\begin{equation}\nonumber
G_{n}^{(stu)} \in \big[-1.96,1.96 \big].
\end{equation}
Each one of these generated 500 coverage probabilities is  based on 500 replications. We then record the proportion of those coverage probabilities  that deviate from the nominal 0.95 by no more than 0.01. This proportion is denoted by $PropG_{n}^{(stu)}(0)$. For the same  generated data, the same proportion, denoted by $PropT_{n}^{(stu)}(0)$,  is also recorded for $T_{n}^{(stu)}(0)$, i.e.,  the classical counterpart of $G_{n}^{(stu)}(0)$.
\par
 The same idea is used to compute the proportions $PropG_{n}^{(stu)}(0)$, $PropG_{n}^{(stu)}(d)$ and $PropG_{n}^{(stu)}(\hat{d})$   and those of their respective classical counterparts $PropT_{n}^{(stu)}(0)$, $PropT_{n}^{(stu)}(d)$ and $PropT_{n}^{(stu)}(\hat{d})$, in Tables 8-12 for  the therein indicated short and long memory processes.

\par
Here again, the choice of $q$ is based on relation (2.14) of  Abadir \emph{et al}. \cite{Abadir  Distaso Giraitis} in each case. More precisely,  we let  $q$ be $ceiling(n^{1/3})$ for the examined  short memory linear processes, and for long memory linear processes with $0< d <0.25$,   we let $q$ be $ceiling(n^{1/(3+4d)})$,  and for $0.25<d<0.5$, we let $q$ be    $ceiling(n^{1/2 -d})$.

\begin{remark}\label{better rate}
It is important to note that, our randomized pivots $G^{{stu}}_n (0)$, $G^{{stu}}_n (d)$ and $G^{{stu}}_n (\hat{d})$, for $\mu=E_X X_1$,  significantly
outperform their  respective classical counterparts $T^{stu}_{n}(0)$, $T^{stu}_{n}(d)$ and $T^{stu}_{n}(\hat{d})$ for short and long memory linear processes. This better performance, most likely, is an indication  that the respective  sampling  distributions of  $G^{{stu}}_n (0)$, $G^{{stu}}_n (d)$ and $G^{{stu}}_n (\hat{d})$ approach that of standard normal at a  faster speed as compared to that  of  their respective  classical counterparts $T^{stu}_{n} (0)$, $T^{stu}_{n} (d)$ and $T^{stu}_{n} (\hat{d})$. In other words, approximating the sampling  distributions of  $G^{{stu}}_n (0)$,  $G^{{stu}}_n (d)$ and $G^{{stu}}_n (\hat{d})$  by  that of standard   normal  most likely  result in   smaller magnitudes  of error in terms the number of observations $n$. The difference in the performance is even more evident when comparing $G_{n}^{{stu}} (d)$ and $G_{n}^{{stu}} (\hat{d})$ to $T^{stu}_{n} (d)$ and  $T^{stu}_{n} (\hat{d})$, respectively, for both Gaussian and non-Gaussian long memory linear processes (cf. Tables 3-6 and Tables 9-12).
\end{remark}

\newpage
\begin{table}
\caption{MA(1): \ $ X_{t}=W_t-0.5\ W_{t-1}$  }
\begin{center}
%Model: $ X_{t}$ is $ARIMA(0,0,0.5) $\\
\begin{tabular}{ c|c|c|c   }
\hline \hline
  Distribution &  $n$ & $PropG_{n}^{{stu}}(0)$ & $PropT_{n}^{stu}(0)$  \\
\hline\hline
% \multirow{3}{*}{$Gaussian\ X_t$} &15  & 0.926 & 0.915 \\
%                                  &20  & 0.936 & 0.933   \\
%                                  &30  & 0.949 & 0.943 \\
%                         \hline \hline
 \multirow{2}{*}{$W_t \ {\substack{d\\=}}\ \ \textmd{Lognormal}(0,1)$} &20&  0.160 & 0.004 \\
                                    &30&  0.210 & 0.156 \\
\hline
\end{tabular}

\end{center}
 \end{table}

\begin{table}
\caption{AR(1):\ $X_t=0.5 \ X_{t-1}+W_t$ }
\begin{center}
%Model: $X_t$ is $ ARIMA(0.5,0,0)$\\
\begin{tabular}{ c|c|c|c   }
\hline \hline
  Distribution &  $n$ &  $PropG_{n}^{{stu}}(0)$ & $PropT_{n}^{stu}(0)$  \\
\hline\hline
% \multirow{3}{*}{$Gaussian\ X_t$} &20  & 0.935 & 0.927 \\
%                                  &30 & 0.947 & 0.941 \\
%                                  &50 & 0.957 & 0.965\\
%                         \hline \hline
 \multirow{2}{*}{$W_t\ {\substack{d\\=}}\ \textmd{Lognormal}(0,1)$} &70& 0.302 & 0.130 \\
                                    &80&  0.536 & 0.384 \\

\hline
\end{tabular}

\end{center}
\end{table}

\begin{table}[H]\label{L1}
\caption{Long Memory with $d=0.2$:\ $X_{t}=(1-B)^{0.2}\ W_t$ }
\begin{center}
%Model: $ X_{t}$ is long memory with $d=0.2$\\
\begin{tabular}{ c|c|c|c   }
\hline \hline
  Distribution &  $n$ & $PropG_{n}^{{stu}}(0.2)$ &  $PropT_{n}^{stu} (0.2)$\\
\hline\hline
% \multirow{2}{*}{$Gaussian\ X_t$} &50  & 0.946 &  0.905  \\
%                                  &70  & 0.953 &  0.928    \\
%                                  %&15  & ? &  ?& \\
%                         \hline \hline
 \multirow{2}{*}{$W_t \ {\substack{d\\=}} \ \textmd{Lognormal}(0,1)$} &150  & 0.686 &0.002  \\
                                                                      &250  & 0.738 & 0.044 \\
\hline
\end{tabular}

\end{center}
\end{table}

\begin{table}[H]\label{Long2}
\caption{Long Memory with $d=0.2$:\ $X_{t}=(1-B)^{0.2} \ W_t$; estimator $\hat{d}$ used}
\begin{center}
\begin{tabular}{ c|c|c|c }
\hline \hline
  Distribution &  $n$ & $PropG_{n}^{{stu}}(\hat{d})$ &  $PropT_{n}^{stu} (\hat{d})$\\
\hline\hline
% \multirow{3}{*}{$Gaussian\ X_t$} &100  & 0.907 &  0.845  \\
%                                  &200  & 0.933 &  0.865    \\
%                                  &300  & 0.955 &  0.908 \\
%                         \hline \hline
 \multirow{2}{*}{$W_t\ {\substack{d\\=}}\ \textmd{Lognormal}(0,1)$} &400  & 0.584 & 0.002  \\
                                    &500  & 0.656 &  0.002  \\
\hline
\end{tabular}

\end{center}
\end{table}

\begin{table}[H]\label{Long3}
\caption{Long Memory with $d=0.4$:\ $X_{t}=(1-B)^{0.4}\ W_t$}
\begin{center}
%Model: $ X_{t}$ is long memory with $d=0.4$\\
\begin{tabular}{ c|c|c|c }
\hline \hline
  Distribution &  $n$ & $PropG_{n}^{{stu}}(0.4)$ &  $PropT_{n}^{stu} (0.4)$\\
\hline\hline

 %\multirow{2}{*}{$Gaussian\ X_t$} &300  & 0.936 &  0.884  \\
 %                                 &400  & 0.947 &  0.901    \\
 %                                 %&300  & 0.955 &  0.908 \\
 %                        \hline \hline
 \multirow{2}{*}{$W_t\ {\substack{d\\=}}\ \textmd{Lognormal}(0,1)$} &300  & 0.214 & 0.000 \\
                                                                    &400  & 0.402 &0.000  \\
\hline
\end{tabular}

\end{center}
\end{table}

\begin{table}[H]\label{G4}
\caption{Long Memory with $d=0.4$:\ $X_{t}=(1-B)^{0.4}\ W_t$; estimator $\hat{d}$ used}
\begin{center}
%Model: $ X_{t}$ is long memory with $d=0.4$\\
\begin{tabular}{ c|c|c|c }
\hline \hline
  Distribution &  $n$ & $PropG_{n}^{{stu}}(\hat{d})$ &  $PropT_{n}^{stu} (\hat{d})$\\
\hline\hline
 \multirow{2}{*}{$W_t\  {\substack{d\\=}}\ \textmd{Lognormal}(0,1)$} &1500  & 0.064 & 0.00 \\
                                                                     &2000  & 0.132 & 0.00 \\
\hline
\end{tabular}

\end{center}
\end{table}

\section{On the bootstrap and  linear processes}\label{Bootstrap and Linear Processes}
In the classical theory of the bootstrap, constructing an asymptotic bootstrap confidence interval for the population mean, based on i.i.d. data,  is done by using the  Student $t$-statistic,  and estimating the underlying percentile of the conditional distribution, given the data, by repeatedly and independently  re-sampling from the  set of data in hand (cf.  Efron and Tibshirani \cite{Efron and Tibshirani}, for example).

\par
  Under certain conditions,    the  cutoff points of the conditional distribution, given the data, of the randomized Student  $t$-statistic, as in (\ref{added 1}) in terms of i.i.d. $X_1,X_2,\ldots$,   are used to
  estimate those  of  the sampling distribution of the  traditional pivot. For more on bootstrapping i.i.d. data we refer to, e.g.,   Davison and Hinkley \cite{Davison  and Hinkley},  Hall \cite{Hall Book} and  Shao and Tu \cite{Shao and Tu}. Considering that the  cutoff points of the randomized $t$-statistic are unknown, they     usually  are  estimated via  drawing $B\geq 2$ independent   bootstrap  sub-samples. The same approach is also taken  when the data form  short or long memory processes.     For  references   on  various bootstrapping methods  to mimic the sampling distribution of  statistics based on  dependent data, and thus, to capture a characteristic of the population,        we refer to H\"{a}rdle \emph{et al}. \cite{Hardle},  Kreiss and Paparoditis \cite{Kreiss and Paparoditis},   Lahiri \cite{Lahiri},  and references therein.

\par
Extending  the i.i.d. based techniques of the bootstrap to fit dependent data by no means can be described as straightforward. Our investigation of the bootstrap, in this section,    sheds  light on some  well known issues that arise  when the bootstrap is applied to long memory  processes.

\par
In this section we   study the problem of bootstrapping linear processes via the same approach that  we used to investigate the asymptotic distribution of $G_{n}$, and its Studentized versions  $G_{n}^{stu}(d)$ and $G_{n}^{stu}(\hat{d})$, the direct randomized pivots for the population mean $\mu$.

\par
In the classical method of conditioning on the data, a bootstrap sub-sample becomes an i.i.d. set of observables   in terms of the classical empirical distribution  even when  the original data are dependent. In comparison,  conditioning on the weights enables us to trace the effect of randomization  by the weights $w^{(n)}_i$ on  the stochastic nature of the original sample.
%The different nature of the two methods to conditioning  becomes  apparent not only in the case of i.i.d. observations,  but also when the data are dependent.

\par
To formally state our results on bootstrapping   linear processes,    we first consider the bootstrapped   sum $ \bar{X}^{*}_n -\bar{X_n}$, where $\bar{X}^{*}_n$ is the mean of a bootstrap sub-sample $X^{*}_1,\ldots,X^{*}_n$ drawn with replacement from the sample $X_i$, $1\leq i \leq n$, of linear processes   as defined in (\ref{eq 1}). Via (\ref{eq 3}), instead of (\ref{eq 4}) as in (\ref{eq 6}), define

\begin{eqnarray}
T_{n}^{*}&:=& \frac{\bar{X}^{*}_n -\bar{X_n}}{\sqrt{\gamma_0 \sum_{j=1}^{n} (\frac{w_{i}^{(n)}}{n}-\frac{1}{n})^2+ 2\sum_{h=1}^{n-1} \gamma_{h} \sum_{j=1}^{n-h} (\frac{w_{j}^{(n)}}{n}-\frac{1}{n})   (\frac{w_{j+h}^{(n)}}{n}-\frac{1}{n}) } }\nonumber\\
&=& \frac{\sum_{i=1}^{n} (\frac{w_{i}^{(n)}}{n}-\frac{1}{n})X_i}{\sqrt{\gamma_0 \sum_{j=1}^{n} (\frac{w_{i}^{(n)}}{n}-\frac{1}{n})^2+ 2\sum_{h=1}^{n-1} \gamma_{h} \sum_{j=1}^{n-h} (\frac{w_{j}^{(n)}}{n}-\frac{1}{n})   (\frac{w_{j+h}^{(n)}}{n}-\frac{1}{n}) } },\nonumber\\
&&\label{eq 5}
\end{eqnarray}
where, as before   in (\ref{eq 3}) and (\ref{eq 4}),  $w_{i}^{(n)}$, $1\leq i \leq n$, are the multinomial random weights of size $n$ with respective probabilities $1/n$, and are independent from the observables in hand.

 \par
 Despite the seeming similarity  of  the bootstrapped statistic $T_{n}^{*}$ to   $G_{n}$ as in (\ref{eq 6}), the two objects are,  in fact,   very different from each other. Apart from the  latter being a direct pivot for the population mean $\mu$, while the former is not, $T^{*}_n$ can only be  used  up to short memory processes. In other words, in case of long memory linear processes,  it fails to converge in distribution to a  non-degenerate limit (cf. Remark \ref{T^* no good long memory}). This is  quite  disappointing when bootstrap is used to capture the sampling distribution of the classical pivot  $T_{n}(d)$ (for $\mu$) of a long memory linear process by that of its bootstrapped  version $T^{*}_n$. It  should also be kept in mind that, in view of Remark \ref{T^* no good long memory},  for  $T_{n}^{*}$  the  natural normalizing sequence, i.e., $Var_{X|w}(\bar{X_{n}^{*}}-\bar{X}_n)=Var_{X|w}\{\sum_{i=1}^n (\frac{w^{(n)}_i}{n}-\frac{1}{n}) X_i\}$ fails to provide the same asymptotic distribution as that of the original statistic $T^{stu}_{n}(d)$, when the data are of long memory.

 \begin{remark}\label{T^* no good long memory}
When dealing with dependent data,    $T^{*}_{n}$ does not preserve the covariance structure of  the original data. This can be explained by observing that the expected values of the   coefficients of the covariance $\gamma_h$, $h\geq 1$, are  $cov_{w}(w^{(n)}_{1} , w^{(n)}_{2} )=-1/n$. As a result (cf. (\ref{eq 4 proofs}) and (\ref{eq 5 proofs}) in the proofs),  as $n\to +\infty$, one has

\begin{equation}\label{eq 2''}
n^{-1} \gamma_{0}  \sum_{j=1}^{n} (w_{j}^{(n)}-1)^2  - \gamma_{0}  =o_{P_{w}}(1),
\end{equation}

\begin{equation}\label{eq 2'}
\sum_{h=1}^{n-1} \gamma_{h} \sum_{j=1}^{n-h} (w_{j}^{(n)}-1)   (w_{j+h}^{(n)}-1)  =o_{P_{w}}(n).
\end{equation}
%The preceding two results, namely (\ref{eq 2''}) and (\ref{eq 2'}), imply that the normalizing sequence of $T^{*}_n$ is $\sqrt{O_{P_{w}}(n)}$,  as $n \to +\infty$.
In view of (\ref{eq 2''}) and (\ref{eq 2'}), for any linear long memory process, $\{X_i, \ i\geq 1\}$, as defined in part (B) of Theorem \ref{CLT G^*^stu},   with a finite and positive variance,  for any $0<d<1/2$, as $n \to +\infty$, we have

\begin{eqnarray*}
&&Var_{X|w} \Big(n^{1/2-d} \big( \bar{X_{n}^{*}}-\bar{X}_n\big)  \Big)\\
&=& Var_{X|w} \Big(n^{1/2-d} \sum_{i=1}^{n} (\frac{w_{i}^{(n)}}{n}-\frac{1}{n})X_i  \Big)  \to 0\ in \ probability-P_w.
\end{eqnarray*}
The latter conclusion        implies  that
$T^{*}_n$   cannot be used for long memory processes.
Hence,   $T^{*}_n$  works  \emph{only} for   short memory  processes.
\end{remark}

\par
In view of   (\ref{eq 2''}) and  (\ref{eq 2'}), for short memory linear processes,    $T^{*}_n$  can, without asymptotic loss of information in probability-$P_w$,   also be defined as
\begin{equation}\label{eqq 5}
T^{*}_{n}:= \frac{\sum_{i=1}^{n} (\frac{w_{i}^{(n)}}{n}-\frac{1}{n})X_i}{\sqrt{\gamma_0 \sum_{j=1}^{n} (\frac{w_{i}^{(n)}}{n}-\frac{1}{n})^2 } }.
\end{equation}
Thus, the two definitions of $T_{n}^{*}$ in  (\ref{eq 5}) and (\ref{eqq 5})  coincide  asymptotically.   We note in passing that the asymptotic equivalence of  (\ref{eq 5}) and (\ref{eqq 5}) does not mean that for a finite number of data, they are equally robust. Obviously, (\ref{eq 5}) is more robust, and it should be used,  for  studying its behavior for   a  finite sample size $n$.
\par
 In the following (\ref{eq 7}), for further study  we present  the Studentized counterpart of  $T^{*}_n$ as defined in (\ref{eqq 5}).
\begin{equation}\label{eq 7}
T_{n}^{*^{stu}}:= \frac{ \sum_{i=1}^{n} (\frac{w_{i}^{(n)}}{n}-\frac{1}{n})X_i}{\sqrt{\bar{\gamma}_0 \sum_{j=1}^{n} (\frac{w_{i}^{(n)}}{n}-\frac{1}{n})^2 } },
\end{equation}

\par
The following two results are respective counterparts of Theorems  \ref{CLT G^*} and \ref{CLT G^*^stu}.

 \begin{thm}\label{CLT T^*}
Suppose that $\{X_i, i \geq 1 \}$ is a stationary linear process as defined  in (\ref{eq 1}) with  $\sum_{k=0}^{\infty} |a_{k}|<+\infty$. Then, as $n\to +\infty$,  we have  for all $t \in \mathbb{R}$,
\begin{equation}\label{eq 9}
P_{X|w} (T_{n}^{*}\leq t) \longrightarrow \Phi(t) \ in \ probability-P_{w}
\end{equation}
and, consequently,
\begin{equation}\label{eq 10}\
P_{X,w} (T_{n}^{*}\leq t) \longrightarrow \Phi(t), \ t \in \mathbb{R}.
\end{equation}

\end{thm}

\begin{thm}\label{CLT T^*^stu}
 Assume that for  the stationary linear process $\{X_i, i \geq 1 \}$, as defined in (\ref{eq 1}),  $\sum_{k=0}^{\infty} |a_k|<+\infty$,
and $E \zeta_{1}^{4}<+\infty$.
 Then, as $n\to +\infty$,  we have  for all $t \in \mathbb{R}$,
\begin{equation}\label{eq 16}
P_{X|w} (T_{n}^{*^{stu}}\leq t) \longrightarrow \Phi(t) \ in \ probability-P_{w}
\end{equation}
and, consequently,
\begin{equation}\label{eq 17}
P_{X,w} (T_{n}^{*^{stu}}\leq t) \longrightarrow \Phi(t).
\end{equation}

\end{thm}

\par
It is also noteworthy that  taking the traditional method of conditioning on the data yields the same conclusion on $T_{n}^{*}$ as the one in Remark \ref{T^* no good long memory} for long memory linear processes. In fact, in case of conditioning on the sample,   recalling that here without loss of generality $\mu=0$, one can see that
\begin{eqnarray*}\nonumber
&&Var(n^{1/2-d}( \bar{X_{n}^*}-\bar{X}_n ) | X_1,\ldots,X_n)\\
&=& Var( n^{1/2-d} \sum_{i=1}^n (\frac{w_{i}^{(n)}}{n}-\frac{1}{n}) X_i | X_1,\ldots,X_n)\\
&=& n^{-2d} (1-\frac{1}{n}) \frac{\sum_{i=1}^n X_{i}^2}{n} - n^{-1-2d} \sum_{h=1}^{n-1} (1-\frac{h}{n}) X_i X_{i+h} \\
&=&o_{P_{X}}(1),\ as \ n\rightarrow +\infty; \ when\ 0< d <1/2.
 \end{eqnarray*}
The preceding convergence to zero takes place when $X_i$'s are of  long memory.
\par
In the literature,  block-bootstrap methods are usually used to modify the randomized $t$-statistic  $T^{*}_n$ so that it should reflect the dependent structure of the data   (cf.  Kreiss and Paparoditis \cite{Kreiss and Paparoditis},   Lahiri \cite{Lahiri} and references therein) that is concealed by the conditional independence of the randomized random variables with common distribution $F_n (x):=n^{-1} \# \{k: 1\leq k \leq n, X_k \leq x \}$, $x \in \mathbb{R}$,  given $X_1,\ldots X_n$.  These  methods are   in contrast to our direct pivot $G_{n}$ as in (\ref{eq 6}), and its Studentized version $G_{n}^{{stu}}(d)$ as defined in (\ref{eq 8}), that can be used \emph{both} for short and long memory processes  without dividing the data into blocks. This is so, since,  the random weights $|w^{(n)}_{i}/n -1/n |$ in $G_{n}$, and in its Studentized version $G^{{stu}}_{n}(d)$ as defined in (\ref{eq 8}), project and preserve
  the covariance structure of  the original sample. To further elaborate on the latter we note that, as $n \to +\infty$, we have
$n E_{w} \big|(w^{(n)}_{1}/n -1/n) (w^{(n)}_{2}/n -1/n) \big| \to 4e^{-2}$ (cf. (\ref{eq 8' proofs}) in the proofs). This, in turn,  implies that, for $0\leq d<1/2$, the   term $n^{1-2d} \sum_{h=1}^{n-1} \gamma_{h} \sum_{j=1}^{n-h} \big|\frac{w_{j}^{(n)}}{n}-\frac{1}{n}\big|   \big| \frac{w_{j+h}^{(n)}}{n}-\frac{1}{n} \big|$ will   be neither zero nor infinity in the limit (cf. (\ref{eq 8'' proofs}) and (\ref{eq 9 proofs}) in the proofs). This means that, unlike $T^{*}_n$,    $G_n$ preserves the covariance structure of the data without more ado.
Hence, $G_{n}$ and its Studentized version $G^{{stu}}_n (d)$, as in (\ref{eq 6}) and  (\ref{eq 8}) respectively,  are natural choices to make inference about the mean of  long memory processes, as well as of  short memory ones.

\subsection*{Advantages of  $G_n$ over bootstrapped  $T_{n}^{*}$  }
In comparing  the use of our direct randomized pivot $G_{n}$ and its Studentized versions  $G^{stu}_n (0)$ and $G_{n}^{stu} (\hat{d})$  %for their respective population means $\mu$
to the $T_{n}^{*}$ based   bootstrap method of constructing confidence intervals, our approach has  a number of  fundamental advantages over the latter. The first advantage is that $G_{n}$ and its Studentized versions can be used without any adjustment of  the  data,  such as dividing them into blocks. The second advantage is that $G_{n}$ and its Studentized versions can be used for both long and short memory linear processes while the bootstrap approach, represented by $T_{n}^{*}$, fails for long memory linear process (cf. Remark\ref{T^* no good long memory}).  The third advantage  concerns the fact that  the Studentized pivots $G^{stu}_n (0)$ and $G_{n}^{stu} (\hat{d})$ are direct  pivots for the population mean  $\mu$. Recall that $T_{n}^{*}$ fails to preserve $\mu$ (cf. Remark \ref{Remark 0}).
The fourth advantage of $G^{stu}_n (0)$  and $G^{stu}_n (\hat{d})$ over the bootstrap methods   is that, apart from the parameter of interest $\mu$, the  other elements of these pivots  are computable based on the  data, and the limiting distribution is standard normal. This is in contrast to using  unknown parameters  as a normalizing sequence  for $ \bar{X^{*}_n}-\bar{X_n}=\sum_{i=1}^n (\frac{w^{(n)}_i}{n}-\frac{1}{n})X_i$ to lift  it to converge to the same limiting distribution, i.e., standard  normal, as that of the original pivots for $\mu$, $T_{n}^{stu}(d)$ and $T_{n}^{stu}(\hat{d})$,   when
block bootstrap is, for example,  applied to  subordinated Gaussian   long memory  processes with Hermit rank one (cf. Theorems 10.2 and 10.3 of  Lahiri \cite{Lahiri}).
Finally, we note that  our approach to making  inference about the population mean $\mu$, based on the pivot $G_n$ and its Studentized versions,  for both short and long memory linear processes,  does  not  require repeated re-sampling from the original data.
This is in contrast to the bootstrap methods, where bootstrap sub-samples are to be drawn repeatedly.

\section{Proofs}\label{Proofs}
In  proving  Theorems   \ref{CLT G^*} and \ref{CLT T^*},   we make use of   Theorem 2.2 of Abadir  \emph{et al}. \cite{Koul} in which the asymptotic normality of sums of deterministically weighted   linear processes  are established. In this context, in Theorems  \ref{CLT G^*} and \ref{CLT T^*},  we view the sums defined in (\ref{eq 3}) and (\ref{eq 4})  as randomly weighted sums of the data on $X$. Conditioning on the weights $w_{i}^{(n)}$s,  we  show that the conditions required for  the deterministic weights in the aforementioned Theorem 2.2 of Abadir  \emph{et al}. \cite{Koul}  hold true in    probability-$P_w$  in this context. The latter, in view of  the characterization of convergence in probability in terms of almost sure convergence of subsequences,   will enable us to conclude  the conditional CLTs, in probability-$P_w$, in Theorems  \ref{CLT G^*} and \ref{CLT T^*}. The unconditional CLTs in terms of the joint distribution $P_{X,w}$ will then  follow from the dominated convergence theorem. Employing Slutsky type arguments, we conclude   Theorems \ref{CLT G^*^stu} and  \ref{CLT T^*^stu} from Theorems \ref{CLT G^*} and  \ref{CLT T^*}, respectively.

\subsection*{Proof of Theorem \ref{CLT G^*}}
In order to prove part (A) of this theorem, we first replace $X_i$ by $X_i/ \gamma_0$, in the definition of the linear process $X_i$ in (\ref{eq 1}), so that, without loss of generality,  we may and shall assume that $Var_{X} (X_i)=1$. By virtue of the latter setup and the assumption that $\gamma_h \geq 0$, for each $ h\geq 1 $, one can readily see that, for each $n\geq 2$,
\begin{eqnarray*}
&&Var_{X|w} \Big( \sum_{j=1}^n  \big| \frac{w_{j}^{(n)}}{n}-\frac{1}{n} \big| (X_j-\mu)  \Big)\\
&=& \sum_{j=1}^n  \big( \frac{w_{j}^{(n)}}{n}-\frac{1}{n} \big)^2 +2 \sum_{h=1}^{n-1} \gamma_h \sum_{j=1}^{n-h} \big|  \big( \frac{w_{j}^{(n)}}{n}-\frac{1}{n} \big) \big( \frac{w_{j+h}^{(n)}}{n}-\frac{1}{n} \big)  \big|\\
&\geq& \sum_{j=1}^n  \big( \frac{w_{j}^{(n)}}{n}-\frac{1}{n} \big)^2.
\end{eqnarray*}
In view of the preceding relation and Theorem 2.2 of  Abadir  \emph{et al}. \cite{Koul}, the proof of part (A) of  Theorem \ref{CLT G^*} follows if we show that, as $n \to +\infty$,
\begin{equation}\label{eq 7 proofs}
\frac{\max_{1\leq i \leq n} \big|   \frac{w_{j}^{(n)}}{n}-\frac{1}{n} \big| } { \sqrt{ \sum_{j=1}^{n} (\frac{w_{i}^{(n)}}{n}-\frac{1}{n})^2+ 2\sum_{h=1}^{n-1} \gamma_{h} \sum_{j=1}^{n-h} \big|\frac{w_{j}^{(n)}}{n}-\frac{1}{n}\big|  \    \big|\frac{w_{j+h}^{(n)}}{n}-\frac{1}{n}\big| }  }=o_{P_{w}}(1).
\end{equation}
The preceding conclusion results from (\ref{eq 2 proofs}), and also  on noting that
\begin{eqnarray*}
&& n \sum_{j=1}^{n} (\frac{w_{i}^{(n)}}{n}-\frac{1}{n})^2+ 2n\sum_{h=1}^{n-1} \gamma_{h} \sum_{j=1}^{n-h} \big|\frac{w_{j}^{(n)}}{n}-\frac{1}{n}\big|  \    \big|\frac{w_{j+h}^{(n)}}{n}-\frac{1}{n}\big|\\
&& \geq n \sum_{j=1}^{n} (\frac{w_{i}^{(n)}}{n}-\frac{1}{n})^2 \to 1\ in\ probability-P_w, \ as \ n \to +\infty. \end{eqnarray*}
The latter result, which  follows from (\ref{eq 4 proofs}), completes  the proof of part (A). $\square$

\par
To establish part (B) of Theorem \ref{CLT G^*}, once again in view of Theorem 2.2 of  Abadir  \emph{et al}. \cite{Koul}, it suffices to show that, as $n \to +\infty$, (\ref{eq 7 proofs}) holds true. Since here we have that $\sum_{h=1}^{\infty} \gamma_h<+\infty$, it suffices to show that, as $n \to +\infty$, we have
%, uniformly in $1\leq h \leq n-1$, as $n \to +\infty$, we have

\begin{equation}\label{eq 8 proofs}
n \sum_{h=1}^n \gamma_h    \ \sum_{j=1}^{n-h} \big|\frac{w_{j}^{(n)}}{n}-\frac{1}{n}\big|  \    \big|\frac{w_{j+h}^{(n)}}{n}-\frac{1}{n}\big| \longrightarrow 4 e^{-2} \sum_{h=1}^{+\infty} \gamma_h \ in \ probability-P_w.
\end{equation}
In order to prove (\ref{eq 8 proofs}), we first note that
\begin{eqnarray}
&&E_{w} \big(  \big|\frac{w_{1}^{(n)}}{n}-\frac{1}{n}\big|      \big|\frac{w_{2}^{(n)}}{n}-\frac{1}{n}\big|   \big)\nonumber\\
&&= -1/n^{3} - 2/n^{2} E_{w}\big\{ (w_{1}^{(n)}-1)(w_{2}^{(n)}-1)   \mathds{1}\big(  (w_{1}^{(n)}-1)(w_{2}^{(n)}-1)<0 \big) \big\}\nonumber\\
&&= -1/n^{3} + 4/n^{2} \big( 1-1/n \big)^n \big(1-1/(n-1) \big)^{n}.\label{eq 8' proofs}
\end{eqnarray}
The preceding relation implies that, as $n \rightarrow +\infty$,
\begin{eqnarray}
&&E_{w}\Big( n \sum_{h=1}^{n-1}\gamma_h  \big(\sum_{j=1}^{n-h} \big|\frac{w_{j}^{(n)}}{n}-\frac{1}{n}\big|  \    \big|\frac{w_{j+h}^{(n)}}{n}-\frac{1}{n}\big| \big) \Big) \nonumber \\
&& = \sum_{h=1}^{n-1}\gamma_h  \frac{-(n-h)}{n^2} + 4 \frac{(n-h)}{n} \big( 1-1/n \big)^n \big(1-1/(n-1) \big)^{n}
 \to 4 \ e^{-2} \sum_{h=1}^{+\infty}\gamma_h . \nonumber\\
 &&\label{eq 8'' proofs}
\end{eqnarray}

\par
By virtue of the preceding result, in order to prove (\ref{eq 8 proofs}), we need to show that, as $n\to +\infty$,
\begin{eqnarray}
&&n \sum_{h=1}^{n-1}\gamma_h \big\{\sum_{j=1}^{n-h} \big|\frac{w_{j}^{(n)}}{n}-\frac{1}{n}\big|    \big|\frac{w_{j+h}^{(n)}}{n}-\frac{1}{n}\big| - ( \frac{-(n-h)}{n^2} +  \frac{4(n-h)}{n} \big( 1-\frac{1}{n} \big)^n \big(1-\frac{1}{n-1} \big)^{n} ) \big\} \nonumber\\
&&=o_{P{w}}(1). \label{eq 9 proofs}
\end{eqnarray}
In order to show the validity of  (\ref{eq 9 proofs}), for  ease of the notation, we define
\begin{equation}\nonumber
b_{n,h}:= \frac{-(n-h)}{n^2} +  \frac{4(n-h)}{n} \big( 1-\frac{1}{n} \big)^n \big(1-\frac{1}{n-1} \big)^{n}
\end{equation}
and, for $\varepsilon>0$, we write
\begin{eqnarray}
&&P_{w} \Big(   \Big| n \sum_{h=1}^{n-1}\gamma_h  \sum_{j=1}^{n-h} \big|\frac{w_{j}^{(n)}}{n}-\frac{1}{n}\big|    \big|\frac{w_{j+h}^{(n)}}{n}-\frac{1}{n}\big|-b_{n,h}  \Big|>\varepsilon    \Big)\nonumber\\
&&\leq \varepsilon^{-1} \sum_{h=1}^{n-1}\gamma_h  E_{w}^{1/2}\Big(  n\sum_{j=1}^{n-h} \big|\frac{w_{j}^{(n)}}{n}-\frac{1}{n}\big|    \big|\frac{w_{j+h}^{(n)}}{n}-\frac{1}{n}\big|-b_{n,h}  \Big)^2.\label{eq 10 proofs}
\end{eqnarray}

\par
We are now to show that $E_{w}\Big(  n\sum_{j=1}^{n-h} \big|\frac{w_{j}^{(n)}}{n}-\frac{1}{n}\big|    \big|\frac{w_{j+h}^{(n)}}{n}-\frac{1}{n}\big|-b_{n,h}  \Big)^2$ approaches zero uniformly in $1\leq h\leq n-1$, as $n\to +\infty$.

\par
Some basic, yet not quite trivial, calculations that also include the use of  the moment generating function of the multinomial distribution show that
\begin{equation}\label{eq 11 proofs}
E_{w}\Big( (\frac{w_{1}^{(n)}}{n}-\frac{1}{n})    (\frac{w_{2}^{(n)}}{n}-\frac{1}{n})  \Big)^2=O(n^{-4})
\end{equation}
and
\begin{eqnarray}
&&E_{w}\Big( \Big|(\frac{w_{1}^{(n)}}{n}-\frac{1}{n})  (\frac{w_{2}^{(n)}}{n}-\frac{1}{n}) (\frac{w_{3}^{(n)}}{n}-\frac{1}{n}) (\frac{w_{4}^{(n)}}{n}-\frac{1}{n}) \Big|  \Big)\nonumber\\
&&= 3/n^{6}- 6/n^{7}+8/n^{4} \big( 1-3/n\big)^{n} \big( 1-1/(n-3) \big)^{n}\nonumber\\
&&+ 8/n^{4} (1-1/n)^{n} \big\{ (n(n-2))\big/(n-1)^{2} -1+ \big( 1-3/(n-1)\big)^{n}  \big\}.\nonumber\\
&&\label{eq 12 proofs}
\end{eqnarray}
By virtue of (\ref{eq 11 proofs}) and (\ref{eq 12 proofs}) we have that

\begin{eqnarray*}
&&E_{w}\Big(  n\sum_{j=1}^{n-h} \big|\frac{w_{j}^{(n)}}{n}-\frac{1}{n}\big|    \big|\frac{w_{j+h}^{(n)}}{n}-\frac{1}{n}\big|-b_{n,h}  \Big)^2 \\
&&= E_{w} \Big(  n\sum_{j=1}^{n-h} \big|\frac{w_{j}^{(n)}}{n}-\frac{1}{n}\big|    \big|\frac{w_{j+h}^{(n)}}{n}-\frac{1}{n}\big|  \Big)^2 - b^{2}_{n,h} \\
&& = n^{2} (n-h) O(n^{-4}) \\
&& + n^2 (n-h) (n-h-1) \Big(3/n^{6}- 6/n^{7}+8/n^{4} \big( 1-3/n\big)^{n} \big( 1-1/(n-3) \big)^{n}\\
&&+ 8/n^{4} (1-1/n)^{n} \big\{ (n(n-2))\big/(n-1)^{2} -1+ \big( 1-3/(n-1)\big)^{n}  \big\}  \Big)-b^{2}_{n,h}.
\end{eqnarray*}
Some further algebra shows that the right hand side of the preceding  relation  can be bounded above by

\begin{equation}\nonumber
n^3 O(n^{-4}) + 3n^{-2} +8/n (1-1/n)^n (1- 1/(n-1))^n \to 0, \ as \ n\to +\infty.
\end{equation}
It is important to note that  the left hand side of the  preceding    convergence   dose not depend on $h$.
 \par
Incorporating now the preceding relation   into (\ref{eq 10 proofs}) yields (\ref{eq 9 proofs}). Now the proof of part (B) of Theorem \ref{CLT G^*} is complete. $\square$
\\
\\
Prior to  establishing the proof of Theorem \ref{CLT G^*^stu}, we first define
\begin{equation}\label{eq 13 proofs}
s_{X}^{2}:=\lim_{n\to +\infty} Var_{X} \big( n^{1/2 -d} \bar{X}_n \big)=\lim_{n\to +\infty} n^{-2d} \big\{ \gamma_0+2\sum_{h=1}^{n-1} \gamma_h (1-h/n)\big\},
\end{equation}
and note  that under regular moment conditions, such as those assumed in Theorem \ref{CLT G^*^stu}, the conclusion  $0<s_{X}^{2}<+\infty$ is valid for $0\leq d <1/2$.

\subsection*{Proof of Theorem \ref{CLT G^*^stu}}
Considering  that in this theorem the data   are linear processes that can be of  short memory, or posses the property of long range dependence, here,   the proofs are given  in a general setup that include both cases.

\par
Prior to stating  the details of the proof of Theorem \ref{CLT G^*^stu}, we note that, when the  $X_i$s  form a long memory process,  in view of the in probability-$P_X$ asymptotic equivalence of the estimator $\hat{d}$ to $d$ as stated in the assumptions of this theorem, we present our proofs for  $G_{n}^{{stu}} (d)$ rather than  for  $G_{n}^{{stu}} (\hat{d})$.

\par
The   proof of both  parts of this theorem will follow if, under their respective conditions,  one shows that for $0\leq d < 1/2$, as  $n,q\to+\infty$ such that $q=O(n^{1/2})$, the following two statements hold true:

\begin{eqnarray}
&&  \frac{n^{1-2d} \gamma_0 \sum_{j=1}^{n} \big( \frac{w_{j}^{(n)}}{n}-\frac{1}{n} \big)^2+2 n^{1-2d} \sum_{h=1}^{n-1} \gamma_h \sum_{j=1}^{n-h} \big|\big( \frac{w_{j}^{(n)}}{n}-\frac{1}{n} \big)\big( \frac{w_{j+h}^{(n)}}{n}-\frac{1}{n} \big) \big| }{n^{-2d}\gamma_0 (1-\frac{1}{n})+ 2 n^{1-2d}  \sum_{h=1}^{n-1} \gamma_h \{ \frac{-(n-h)}{n^{3}}+\frac{4(n-h)}{n^2} (1-\frac{1}{n})^n (1-\frac{1}{n-1})^n \}     } \nonumber\\
&& \longrightarrow 1\ in \ probability-P_{w},  \label{eq 29 proofs}\\
&&P_{X|w} \Big\{ \Big|  \frac{\frac{n \bar{\gamma}_0 \sum_{j=1}^{n} \big( \frac{w_{j}^{(n)}}{n}-\frac{1}{n} \big)^2}{q^{2d}}+2 q^{-1-2d} \sum_{h=1}^{q} \bar{\gamma}_h \sum_{j=1}^{q-h} \big| \big( {w_{j}^{(n)}}-{1} \big)\big( {w_{j+h}^{(n)}}-{1} \big)\big| } { \frac{\bar{\gamma}_0 (1-1/n)}{q^{2d}}+ 2 q^{-2d} \sum_{h=1}^{n-1}  \bar{\gamma}_h (1-\frac{h}{q}) \{ \frac{-1}{n} + 4 (1-\frac{1}{n})^n (1- \frac{1}{n-1})^n \}  }-1   \Big|>\varepsilon  \Big\}\nonumber\\
&&=o_{P_{w}}(1),\label{eq 30 proofs}
\end{eqnarray}
where, $\varepsilon$ is an arbitrary positive number.
\par
Due to the following two conclusions,  namely (\ref{eq 31 proofs}) and (\ref{eq 32 proofs}),   (\ref{eq 29 proofs}) and (\ref{eq 30 proofs}) will, in turn,  imply Theorem \ref{CLT G^*^stu}. We have, as $n\rightarrow +\infty$,

\begin{equation}\label{eq 31 proofs}
n^{-2d}\big\{{\gamma}_0+2\sum_{h=1}^{n-1} {\gamma}_h (1-h/n)\big\}\to s_{X}^{2} \ in \ probability-P_{X},
\end{equation}
and, as $n,q\to +\infty$ such that $q=O(n^{1/2})$,
\begin{equation}\label{eq 32 proofs}
q^{-2d}\big\{\bar{\gamma}_0+2\sum_{h=1}^{q} \bar{\gamma}_h (1-h/q)\big\} \to s_{X}^{2} \ in \ probability-P_{X},
\end{equation}
where, $s_{X}^{2}$ is as defined in (\ref{eq 13 proofs}). In the context of Theorem \ref{CLT G^*^stu}, the   conclusion    (\ref{eq 32 proofs}) results from   Theorem 3.1  of Giraitis  \emph{et al}. \cite{Giraitis et al.}. This is so, since,  in Theorem \ref{CLT G^*^stu} we  assume  that the data have a finite forth moment,  $n,q \to +\infty$ in such a way that $q=O(n^{1/2})$, and   in the case of long memory, in part (B) of  Theorem \ref{CLT G^*^stu}  we consider long memory linear processes for which we have   $a_i\sim c i^{d-1}$, as $i \rightarrow +\infty$.

\par
In order to prove (\ref{eq 29 proofs}), we note that, as  $n \to +\infty$,
\begin{eqnarray*}\nonumber
&&n^{-2d} \gamma_0 (1-\frac{1}{n})+ 2 n^{1-2d}  \sum_{h=1}^{n-1} \gamma_h \{ \frac{-(n-h)}{n^{3}}+\frac{4(n-h)}{n^2} (1-\frac{1}{n})^n (1-\frac{1}{n-1})^n \} \\
&&\longrightarrow \left\{
                                                                                    \begin{array}{ll}
                                                                                      (1-4 e^{-2})\gamma_0+ 4 e^{-2} s_{X}^{2}, & \hbox{when\ $d=0$;} \\
                                                                                      4 e^{-2} s_{X}^{2}, & \hbox{when\ $0<d<1/2$.}
                                                                                    \end{array}
                                                                                  \right.
\end{eqnarray*}
Considering that here we have  $\limsup_{n \to +\infty}  n^{-2d} \sum_{h=1}^{n} |\gamma_{h}|<+\infty$,     (\ref{eq 4 proofs}) and (\ref{eq 9 proofs}) imply (\ref{eq 29 proofs}), as $n \to +\infty$.

\par
In order to establish (\ref{eq 30 proofs}), in view of (\ref{eq 4 proofs}) and the  fact that, under the conditions of Theorem \ref{CLT G^*^stu}, as $n\to +\infty$, $\bar{\gamma}_{0}- \gamma_0=o_{P_{X}}(1)$,  we conclude that, as $n \to +\infty$,

\begin{equation*}
P_{X|w} \big( \big| n^{1-2d} (\frac{q}{n})^{-2d} \bar{\gamma_0} \sum_{i=1}^n \big( \frac{w_{i}^{(n)}}{n}-\frac{1}{n} \big)^2- q^{-2d} \bar{\gamma}_{0} (1-1/n) \big|> \varepsilon   \big) \to 0\ in \ probability-P_w,
\end{equation*}
where $0\leq d <1/2$ and   $\varepsilon>0$ is arbitrary.

\par
We proceed with the proof of (\ref{eq 30 proofs}) by showing that, as  $n,q\to +\infty$, the following relation holds true: for arbitrary $\varepsilon_1,\varepsilon_2>0$,  as $n,q \to +\infty$, in such a way that $q=O(n^{1/2})$,
\begin{equation}\label{needed 1}
P_{w}\big\{ P_{X|w} (q^{-2d} \big| \sum_{h=1}^q    \bar{\gamma_{h}} B_{n,q}\big|> 2\varepsilon_1  )>2\varepsilon_2  \big\} \rightarrow 0,
\end{equation}
where

\begin{eqnarray*}\label{eq 35 proofs}
B_{n,q}(h)&:=& q^{-1} \sum_{j=1}^{q-h} \big| \big( w_{j}^{(n)}-1 \big)\big( w_{j+h}^{(n)}-1 \big)\big|-b_{n,q,h},\\
b_{n,q,h}&:=&E_{w} \big( q^{-1} \sum_{j=1}^{q-h} \big| \big( w_{j}^{(n)}-1 \big)\big( w_{j+h}^{(n)}-1 \big)\big|  \big)  \\
&=& \frac{-(q-h)}{n q} +  \frac{4(q-h)}{q} \big( 1-\frac{1}{n} \big)^n \big(1-\frac{1}{n-1} \big)^{n}.
\end{eqnarray*}

\par
In order to establish (\ref{needed 1}), without loss of generality, we first assume $\mu=E_{X} X_{1}=0$,  and for each $1\leq h\leq q$ define
\begin{equation}\label{def. gamma*}
\gamma_{h}^*:=\frac{1}{n} \sum_{i=1}^{n-h}  X_{i} X_{i+h}.
\end{equation}
\par
Observe  now that the left hand side of (\ref{needed 1}) is bounded above by
\begin{eqnarray}
&&P_{w}\big\{ P_{X|w}\big( q^{-2d}\big|  \sum_{h=1}^q (\bar{\gamma}_{h}-\gamma_{h}^*) B_{n,q}(h)   \big|>\varepsilon_1   \big)>\varepsilon_2 \big\}\nonumber\\
&+&  P_{w}\big\{ P_{X|w}\big( q^{-2d} \big|  \sum_{h=1}^q \gamma_{h}^* B_{n,q}(h)   \big|>\varepsilon_1   \big)>\varepsilon_2 \big\}.\label{eq 21 proofs}
\end{eqnarray}
We now show that the first term in (\ref{eq 21 proofs}), i.e., the remainder, is asymptotically negligible. To do so, we  note that we have
\begin{eqnarray}
\sum_{h=1}^q (\bar{\gamma}_{h}-\gamma_{h}^*) B_{n,q}(h)&=&-\frac{\bar{X}_n}{n} \sum_{h=1}^q  B_{n,q}(h)  \sum_{i=1}^{n-h} X_i - \frac{\bar{X}_n}{n} \sum_{h=1}^q  B_{n,q}(h) \sum_{i=1}^{n-h} X_{i+h}\nonumber\\
&&+
\bar{X}^2 \sum_{h=1}^q B_{n,q}(h)\nonumber\\
&\sim& - \bar{X}^2  \sum_{h=1}^q B_{n,q}(h)\ uniformly \ in \ h\ in \ probability-P_{X|w},\nonumber\\
&&\label{eq 22 proofs}
\end{eqnarray}
 where, in the preceding conclusion generically,  $Y_{n}\sim Z_{n}$ in probability-$P$ stands for the in probability-$P$ asymptotic equivalence of the sequences of  random variables $Y_n$ and $Z_n$. The approximation in (\ref{eq 22 proofs}) is true since, for example, for $\varepsilon>0$
\begin{eqnarray*}
P_{X}\big( \cup_{1\leq h \leq q} \big|  \bar{X}_{n} - \frac{\sum_{i=1}^{n-h} X_i}{n} \big|>\varepsilon \big)&\leq&q  P_{X}\big(  \big|  \frac{\sum_{i=n-h+1}^{n} X_i}{n} \big|>\varepsilon \big)\label{newly added 1}\\
&\leq& \varepsilon^{-4}  q \frac{(h-1)^4}{n^4} E(X^{4}_1)\nonumber\\
&\leq&\varepsilon^{-4} \frac{q^5}{n^4} E_{X}(X^{4}_1)\rightarrow 0,\ as \ n\to +\infty.\nonumber
\end{eqnarray*}
The preceding is true since $1\leq h \leq q$ and $q=O(n^{1/2})$, as $n,q \to +\infty$.
\par
We note that   for $0\leq d <1/2$, as $n\rightarrow +\infty$, we have that  $n^{1/2-d} \bar{X}_{n}=O_{P_{X}}(1)$. The latter conclusion, in view of  the equivalence in (\ref{eq 22 proofs}), implies that, for each $\varepsilon_1, \varepsilon_2>0$, there exists   $\varepsilon>0$ such that
\begin{eqnarray}
 &&P_{w}\big\{ P_{X|w}\big( q^{-2d} \big|  \sum_{h=1}^q (\bar{\gamma}_{h}-\gamma_{h}^*) B_{n,q}(h)   \big|>\varepsilon_1   \big)>\varepsilon_2 \big\}\nonumber\\
&\sim& P_{w} \big\{ \frac{q^{-2d}}{n^{1-2d}} \sum_{h=1}^q \big|B_{n,q}(h)\big|>\varepsilon\big\}\nonumber \\
&\leq& \varepsilon^{-1} \frac{q^{-2d}}{n^{1-2d}} \sum_{h=1}^q E_{w}\big(\big|B_{n,q}(h)\big|\big). \label{change 1}
\end{eqnarray}
Observing now  that  $\sup_{n\geq 2}\sup_{1\leq h\leq q} E_{w}\big(\big|B_{n,q}(h)\big|\big)\leq 10$,  we can bound the preceding relation above by
\begin{equation*}
10 \ \varepsilon^{-1}  \frac{q^{1-2d}}{n^{1-2d}}\longrightarrow 0,
\end{equation*}
as $n,q \to +\infty$ in such away that $q=O(n^{1/2})$. This means that the first term in (\ref{eq 21 proofs}) is asymptotically negligible and, as a result,  (\ref{needed 1}) follows when  the second term in the former  relation is also asymptotically negligible.  To prove this negligibility, we  first define
\begin{equation}\label{def. of gamma**}
\gamma_{h}^{**}:=\frac{1}{n} \sum_{i=1}^n X_i X_{i+h}.
\end{equation}
Now, observe that
\begin{eqnarray}
P_{X}\big\{  \cup_{1\leq h \leq q} |\gamma_{h}^{**}-\gamma_{h}^{*}|>\varepsilon \big\}&\leq& q P\big\{   \frac{1}{n}| \sum_{i=n-h+1}^n X_i X_{i+h}|>\varepsilon \big\}\nonumber\\
&\leq& \varepsilon^{-2} \frac{q^3}{n^2} E_{X} (X^{4}_1)\to 0,\nonumber
\end{eqnarray}
as $n,q\to +\infty$ such that $q=O(n^{1/2})$,  hence, as $n,q\to +\infty$ such that $q=O(n^{1/2})$, using a similar argument to arguing   (\ref{eq 21 proofs}) and (\ref{change 1}),  with $\gamma^{*}_h$ replacing $\bar{\gamma}_h$ and $\gamma^{**}_h$ replacing $\gamma^{*}_h$ therein,  we  arrive at
\begin{eqnarray*}
&&P_{w}\big\{ P_{X|w}\big( q^{-2d}\big|  \sum_{h=1}^q \gamma_{h}^{*} B_{n,q}(h)   \big|>\varepsilon_1   \big)>\varepsilon_2 \big\}\\
&\sim& P_{w}\big\{ P_{X|w}\big( q^{-2d}\big|  \sum_{h=1}^q \gamma_{h}^{**} B_{n,q}(h)   \big|>\varepsilon_1   \big)>\varepsilon_2 \big\}.
\end{eqnarray*}

\par
Therefore, in order to prove (\ref{needed 1}), it suffices to show that, as $n,q\to +\infty$ so that $q=O(n^{1/2})$,
\begin{equation*}\label{eq 36 proofs}
P_{w}\big\{ P_{X|w}\big( q^{-2d}\big|  \sum_{h=1}^q \gamma_{h}^{**} B_{n,q}(h)   \big|>\varepsilon_1   \big)>\varepsilon_2 \big\}\to 0,
\end{equation*}
where  $\gamma_{h}^{**}$ is defined in (\ref{def. of gamma**}). The  latter relation, in turn,   follows from the following two conclusions: as $n,q \to +\infty$ so that $q=O(n^{1/2})$,

\begin{equation}\label{eq 37 proofs}
\sup_{1\leq h,h^\prime \leq q}  E_{w}\big( \big|B_{n,q}(h) B_{n,q}(h^{\prime})  \big| \big)=o(1)
\end{equation}
and
\begin{equation}\label{eq 38 proofs}
q^{-4d} \sum_{h=1}^q \sum_{h^{\prime}=1}^q  \big| E_{X}(\gamma_{h}^{**}  \gamma_{h^{\prime}}^{**}) \big|=O(1).
\end{equation}

\par
To prove (\ref{eq 37 proofs}),  we use the Cauchy inequality to write

\begin{eqnarray*}
&&E_{w}\big( \big|B_{n,q}(h) B_{n,q}(h^{\prime})  \big| \big)\\
&\leq& E_{w}\big( B_{n,q}(h) \big)^2\\
&\leq&   \frac{q-h}{q^2} E_{w}\big( (w^{(n)}_{1}-1) (w^{(n)}_{2}-1)   \big)^2\\
&+& \frac{(q-h)(q-h-1)}{q^2} E_{w} \big|   (w^{(n)}_{1}-1) (w^{(n)}_{2}-1) (w^{(n)}_{3}-1) (w^{(n)}_{4}-1)    \big| \\
&-& b_{n,q,h}^{2}\\
&\leq& (q-1)/q^{2} O(1)+ 3/n^2 +8/n (1-1/n)^n (1-1/(n-1))^n.
%&=&  \frac{q-h}{q^2} O(1)+ 3 \frac{(q-h)(q-h-1)}{q^2 n^2} - 6 \frac{(q-h)(q-h-1)}{q^2 n^3} \\
%&+& 8 \frac{(q-h)(q-h-1)}{q^2} (1-\frac{3}{n})^n (1- \frac{1}{n-3})^n\\
%&+& 8 \frac{(q-h)(q-h-1)}{q^2} (1-\frac{1}{n})^n \big\{\frac{ n(n-2)}{(n-1)^2}-1 +(1-\frac{3}{n-1})^n \big\}\\
%&-& b_{n,q,h}^{2}.
\end{eqnarray*}

\par
We note that the right hand side of the preceding relation does not depend on $h$ and it approaches zero as $n\to +\infty$. The latter conclusion implies (\ref{eq 37 proofs}).

\par
In order to establish (\ref{eq 38 proofs}),  we define
 \begin{equation*}
H:=\lim_{s\to +\infty} s^{-2d} \sum_{\ell =-s}^{s} |\gamma_{\ell}|.
\end{equation*}
Observe that $H<+\infty$.
We now carry on with the proof of (\ref{eq 38 proofs}) using  a generalization of an argument used in the proof of Proposition 7.3.1 of Brockwell  and Davis \cite{Brockwell and Davis} as follows:

\begin{eqnarray}
&&q^{-4d}\sum_{h=1}^q  \sum_{h^\prime=1}^q \big|  E_{X}(\gamma_{h}^{**}  \gamma_{h^{\prime}}^{**}) \big|\nonumber\\
 &\leq& q^{-2d}\sum_{h=1}^q  \big|\gamma_h \big|\ q^{-2d} \sum_{h^\prime=1}^q \big|\gamma_{h^\prime}|\nonumber\\
&+&(\frac{q}{n})^{1-2d} n^{-2d} \sum_{k=-n}^n \big| \gamma_h \big|\ q^{-2d} \sum_{L=-q}^q  \big| \gamma_{k+L} \big|\nonumber \\
&+& \frac{1}{n} \sum_{k=-n}^n \ q^{-2d} \ \sum_{h^\prime=1}^q  \big| \gamma_{k+h^\prime} \big| \ q^{-2d} \ \sum_{h=1}^q \big| \gamma_{k-h} \big|\nonumber\\
&+& \frac{q^{-2d}}{n^{1-2d}}  n^{-2d} \sum_{i=1}^n \sum_{k=-n}^n a_i a_{i+k}\  q^{-d} \sum_{h=1}^q a_{i+h} \ q^{-d} \sum_{h^\prime=1}^q a_{i+k-h^\prime}.\label{eq 39 prroofs}
\end{eqnarray}
It is easy to see that, as $n\to +\infty$, and consequently $q\to +\infty$, the right hand side of the inequality  (\ref{eq 39 prroofs}) converges to the finite limit $3 H^2$. Now the proof of (\ref{eq 38 proofs}) and also that of Theorem \ref{CLT G^*^stu} are complete. $\square$

\subsection*{Proof of Corollary \ref{corollary 1}}
Due to the similarity of parts (A) and (B), we only give the proof for part  (A) of  Corollary \ref{corollary 1}.
\par
In order to establish part (A), we first construct an asymptotic  $1-\alpha$ size one-sided randomized  confidence bound for the parameter $\mu_{X}=E_{X} X$
using part (B) of Theorem \ref{CLT G^*^stu}, as follows:
\begin{equation}\label{eq 40 prroofs}
 \mu_{X}\geq \frac{ \sum_{i=1}^n |\frac{w_{i}^{(n)}}{n}-\frac{1}{n}|X_i -  D^{1/2}_{n,q,\hat{d}} z_{1-\alpha} }{ \sum_{j=1}^n |\frac{w_{j}^{(n)}}{n}-\frac{1}{n}|}.
\end{equation}
Now, since the function $\mathcal{G}$ is an increasing function, we conclude that (\ref{eq 40 prroofs}) is equivalent to having
\begin{equation}\nonumber
\mathcal{G}(\mu_{X})\geq \mathcal{G}(\frac{ \sum_{i=1}^n |\frac{w_{i}^{(n)}}{n}-\frac{1}{n}|X_i - D^{1/2}_{n,q,\hat{d}} z_{1-\alpha} }{ \sum_{j=1}^n |\frac{w_{j}^{(n)}}{n}-\frac{1}{n}|}).
\end{equation}
Employing   Jenssen's inequality at this stage  yields   conclusion   (A) of Corollary \ref{corollary 1}.
Now the proof of Corollary \ref{corollary 1} is complete. $\square$

\subsection*{Proof of Theorem \ref{CLT T^*}}
Without loss of generality here, we assume that $\mu=0$, and note that
\begin{eqnarray}
&&Var_{X|w}\big( \sum_{i=1}^{n} (\frac{w_{i}^{(n)}}{n}-\frac{1}{n}) X_i \big)\nonumber\\
&&= \gamma_0 \sum_{j=1}^{n} (\frac{w_{j}^{(n)}}{n}-\frac{1}{n})^2+ 2\sum_{h=1}^{n-1} \gamma_{h} \sum_{j=1}^{n-h} (\frac{w_{j}^{(n)}}{n}-\frac{1}{n})   (\frac{w_{j+h}^{(n)}}{n}-\frac{1}{n}).\nonumber
\end{eqnarray}
Now, in view of Theorem 2.2 of Abadir  \emph{et al}. \cite{Koul}, it suffices to show that, as $n \to +\infty$,
\begin{equation}\label{eq 1 proofs}
\frac{\max_{1\leq i \leq n} \big(\frac{w_{i}^{(n)}}{n}-\frac{1}{n}  \big)^2 }{\gamma_0 \sum_{j=1}^{n} (\frac{w_{j}^{(n)}}{n}-\frac{1}{n})^2+ 2\sum_{h=1}^{n-1} \gamma_{h} \sum_{j=1}^{n-h} (\frac{w_{j}^{(n)}}{n}-\frac{1}{n})   (\frac{w_{j+h}^{(n)}}{n}-\frac{1}{n})}=o_{P_{w}}(1).
\end{equation}
Noting  that $\gamma_{0} >0$, the proof of the preceding statement results from the following two conclusions: as $n\to +\infty$,
\begin{equation}\label{eq 2 proofs}
 n  \max_{1\leq i \leq n} \big(\frac{w_{i}^{(n)}}{n}-\frac{1}{n}  \big)^2= o_{P_{w}}(1),
\end{equation}
and
\begin{equation}\label{eq 3 proofs}
n \gamma_0 \sum_{j=1}^{n} (\frac{w_{j}^{(n)}}{n}-\frac{1}{n})^2+ 2n\sum_{h=1}^{n-1} \gamma_{h} \sum_{j=1}^{n-h} (\frac{w_{j}^{(n)}}{n}-\frac{1}{n})   (\frac{w_{j+h}^{(n)}}{n}-\frac{1}{n})- \gamma_{0}= o_{P_{w}}(1).
\end{equation}

\par
To prove (\ref{eq 2 proofs}), for $\varepsilon>0$,  in what follows we  employ Bernstien's inequality and  write
\begin{eqnarray*}
P_{w}\big(  \max_{1\leq i \leq n} \big|\frac{w_{i}^{(n)}}{n}-\frac{1}{n}  \big|>\frac{\varepsilon}{\sqrt{n} }\big)&\leq& n  P_{w}\big(   \big|\frac{w_{1}^{(n)}}{n}-\frac{1}{n}  \big|>\frac{\varepsilon}{\sqrt{n} }\big)\\
&\leq& n \exp\{  - n^{1/2} \frac{\varepsilon^2}{n^{-1/2}+\varepsilon} \}=o(1),
\end{eqnarray*}
as $n\to +\infty$. Now the proof of (\ref{eq 2 proofs}) is complete.

\par
Considering that  here we have $\sum_{h=1}^{\infty} \gamma_h<+\infty$, the proof  of (\ref{eq 3 proofs}) will follow from the following two statements: as $n\to +\infty$,
\begin{equation}\label{eq 4 proofs}
n \sum_{i=1}^{n} \big(\frac{w_{i}^{(n)}}{n}-\frac{1}{n}  \big)^2-(1-1/n) = o_{P_{w}}(1)
\end{equation}
and
\begin{equation}\label{eq 5 proofs}
n\sum_{h=1}^{n-1} \gamma_{h}  \sum_{j=1}^{n-h} \big(\frac{w_{j}^{(n)}}{n}-\frac{1}{n}  \big) \big(\frac{w_{j+h}^{(n)}}{n}-\frac{1}{n}  \big)=o_{P_{w}}(1).
\end{equation}
To prove (\ref{eq 4 proofs}), with  $\varepsilon>0$, we first use Chebyshev's inequality followed by some algebra involving the use of the moment generating function of the  multinomial distribution to write
\begin{eqnarray*}
&&P_{w} \big(  \Big| n \sum_{i=1}^{n} \big( \frac{w^{(n)}_i}{n}-\frac{1}{n}  \big)^2 - (1- 1/n)  \Big|>\varepsilon   \big)\\
&\leq&   \varepsilon^{-2} n^2  E_{w}\big( \sum_{i=1}^{n} (\frac{w^{(n)}_{i}}{n}-\frac{1}{n})^2 - \frac{(1-\frac{1}{n})}{n}   \big)^2\\
&\leq& \varepsilon^{-2} n^2  (1-\frac{1}{n})^{-2}  \Big\{ \frac{1-\frac{1}{n}}{n^6  } + \frac{(1-\frac{1}{n})^4}{n^3}  + \frac{(n-1)(1-\frac{1}{n})^2}{n^4 } +  \frac{4(n-1)}{n^4}       +\frac{1}{n^2}\nonumber\\
&& - \frac{1}{n^3} + \frac{n-1}{n^6} + \frac{4(n-1)}{n^5} - \frac{(1-\frac{1}{n})^2}{n^2} \Big\}
\longrightarrow 0,\ \textrm{as}\ n\to +\infty.
\end{eqnarray*}
The latter completes the proof of (\ref{eq 4 proofs}).

\par
In order to establish (\ref{eq 5 proofs}), with  $\varepsilon>0$, we  write
\begin{eqnarray}
&&P_{w}\big( n \big|\sum_{h=1}^n \gamma_{h} \sum_{j=1}^{n-h} \big(\frac{w_{j}^{(n)}}{n}-\frac{1}{n}  \big) \big(\frac{w_{j+h}^{(n)}}{n}-\frac{1}{n}   \big) \big| >\varepsilon \big)\nonumber\\
&\leq& \varepsilon^{-1} \sum_{h=1}^{n-1} |\gamma_{h}| E^{1/2}_{w}\Big( n\sum_{j=1}^{n-h}  \big( \frac{w_{j}^{(n)}}{n}-\frac{1}{n}  \big) \big(\frac{w_{j+h}^{(n)}}{n}-\frac{1}{n} \big)  \Big)^2. \label{eq 6 proofs}
\end{eqnarray}
Observe now that

\begin{eqnarray*}
&&E_{w}\Big( n\sum_{j=1}^{n-h}  \big( \frac{w_{1}^{(n)}}{n}-\frac{1}{n}  \big) \big(\frac{w_{2}^{(n)}}{n}-\frac{1}{n} \big)  \Big)^2     \\
&=& n^2 (n-h) E_{w} \Big(  \big( \frac{w_{1}^{(n)}}{n}-\frac{1}{n}  \big) \big(\frac{w_{2}^{(n)}}{n}-\frac{1}{n} \big)   \Big)^2   \\
&+&n^2 (n-h)(n-h-1) E_{w} \Big(  \big( \frac{w_{1}^{(n)}}{n}-\frac{1}{n}  \big) \big(\frac{w_{2}^{(n)}}{n}-\frac{1}{n} \big)\big( \frac{w_{3}^{(n)}}{n}-\frac{1}{n}  \big) \big(\frac{w_{4}^{(n)}}{n}-\frac{1}{n} \big)     \Big)\\
&\leq& n^3 O(n^{-4}) +n^3 O(n^{-6})\to 0,\ as\ n\to +\infty.
\end{eqnarray*}
The preceding conclusion is true, since  $ E_{w} \Big(  \big( \frac{w_{1}^{(n)}}{n}-\frac{1}{n}  \big) \big(\frac{w_{2}^{(n)}}{n}-\frac{1}{n} \big)  \Big)^{2} =O(n^{-4})$ and
$ E_{w} \Big(  \big( \frac{w_{1}^{(n)}}{n}-\frac{1}{n}  \big) \big(\frac{w_{2}^{(n)}}{n}-\frac{1}{n} \big) \big( \frac{w_{3}^{(n)}}{n}-\frac{1}{n}  \big)  \big( \frac{w_{4}^{(n)}}{n}-\frac{1}{n}  \big) \Big)
= O({n^{-6}})$. Incorporating now the latter two results into (\ref{eq 6 proofs}), the conclusion (\ref{eq 5 proofs}) follows.  Now  the proof of Theorem \ref{CLT T^*} is complete. $\square$

\subsection*{Proof of Theorem \ref{CLT T^*^stu}}
In order   to prove   Theorem \ref{CLT T^*^stu}, using a Slutsky type argument,     it suffices  to show that the   Studentizing sequence of $T_{n}^{*^{stu}}$,  asymptotically in $n$,  in a hierarchical   way,   coincides with the right normalizing sequence, i.e., with the one in the denominator  of $T_{n}^{*}$ defined in (\ref{eqq 5}).
\par
Considering that, as $n \to +\infty$, we have that  $\bar{\gamma}_0 -\gamma_{0}=o_{P_{X}}(1)$, where $0< \gamma_0 <+\infty$, the proof of this theorem follows if, for $\varepsilon_1, \varepsilon_2 >0$, we show that
\begin{equation*}
P_{w} \big\{ P_{X|w} \big( \bar{\gamma}_0 | n\sum_{i=1}^n  (\frac{w_{i}^{(n)}}{n} -\frac{1}{n}  )-(1-1/n)|>\varepsilon_1     \big)>\varepsilon_2  \big\}=o(1), \ as\ n\rightarrow +\infty.
\end{equation*}
\par
To establish the preceding relation, we note that its left hand side is bounded above by
\begin{eqnarray*}
&&P_{w}\big\{ E_{X} (\bar{\gamma}_0) \big( |n \sum_{i=1}^n  (\frac{w_{i}^{(n)}}{n} -\frac{1}{n}  )-
(1-1/n)|     \big)>\varepsilon_1 \varepsilon_2  \big\}\\
&& \leq P_{w}\big\{  \big( |n \sum_{i=1}^n  (\frac{w_{i}^{(n)}}{n} -\frac{1}{n}  )-(1-1/n)|     \big)>\frac{\varepsilon_1 \varepsilon_2}{\gamma_0}  \big\}.
\end{eqnarray*}
The rest of the proof is similar to that of (\ref{eq 4 proofs}).
Now the proof of this theorem is complete. $\square$

% d=0.4 with n=1500 and dhat 0.064 AND 0.00
% d=0.4 with n=2000 and dhat 0.132 AND 0.00
\end{document}